\documentclass{amsart}
\usepackage{amssymb}
\usepackage{upref}

\newtheorem{theorem}{Theorem}[section]
\newtheorem{lemma}[theorem]{Lemma}
\newtheorem{proposition}[theorem]{Proposition}
\newtheorem{corollary}[theorem]{Corollary}
\newtheorem{conjecture}[theorem]{Conjecture}
\newtheorem{question}[theorem]{Question}
\newtheorem{problem}[theorem]{Problem}
\theoremstyle{definition}
\newtheorem{definition}[theorem]{Definition}

\newtheorem{example}[theorem]{Example}

\DeclareMathOperator{\loc}{loc}
\DeclareMathOperator{\Hom}{Hom}
\DeclareMathOperator{\lub}{lub}
\renewcommand{\lub}{\bigvee}

\newcommand{\bousclass}[1]{\langle #1 \rangle}
\newcommand{\bcdual}[1]{I \mspace{-2mu} {#1}}
\newcommand{\doublebcdual}[1]{I^{2} \mspace{-2mu} {#1}}
\newcommand{\BP}{B\!P}
\newcommand{\bpj}{B\!P\!J}
\newcommand{\HFp}{H \mspace{-2mu} \zp}
\newcommand{\HZ}{H \mspace{-2mu} \mathbf{Z}}
\newcommand{\MSp}{M\!Sp}
\newcommand{\DL}{\mathbf{DL}}
\newcommand{\BA}{\mathbf{BA}}
\newcommand{\cBA}{\mathbf{cBA}}
\newcommand{\lattice}{\mathbf{B}}
\newcommand{\locallattice}{L_{D}\lattice}
\newcommand{\strange}{J}
\newcommand{\swdual}[1]{D \mspace{-2mu} {#1}}
\newcommand{\Smash}{\wedge}

\newcommand{\zp}{\mathbf{F}_{p}}

\newcommand{\cat}[1]{\mathcal{#1}}
\newcommand{\mathcolon}{\colon\,}
\newcommand{\suchthat}{\,\vert\,}
\newcommand{\card}[1]{\beth_{#1}}
\newcommand{\meet}{\curlywedge}

\newcommand{\ulp}{\textup{(}}
\newcommand{\urp}{\textup{)}}
\newcommand{\uc}{\textup{:}}
\newcommand{\usc}{\textup{;}}

\newcommand{\Strut}[1]{\rule{0in}{#1}}
\newcommand{\lbp}{\left( \Strut{0.125in}}
\newcommand{\rbp}{\right)}

\begin{document}
 
\title{The structure of the Bousfield lattice}


\date{\today}

\author{Mark Hovey}
\address{Department of Mathematics \\ Wesleyan University
\\ Middletown, CT 06459}
\email{hovey@member.ams.org}
\thanks{Research partially supported by a National Science Foundation
grant}

\author{John H. Palmieri}
\address{Department of Mathematics \\ University of Notre Dame \\
Notre Dame, IN 46556}
\email{palmieri@member.ams.org}
\thanks{Research partially supported by National Science Foundation
grant DMS-9407459}

\date{\today}

\begin{abstract}
Using Ohkawa's theorem that the collection $\lattice$ of Bousfield
classes is a set, we perform a number of constructions with Bousfield
classes.  In particular, we describe a greatest lower bound operator;
we also note that a certain subset $\DL$ of $\lattice$ is a frame, and
we examine some consequences of this observation.  We make several
conjectures about the structure of $\lattice$ and $\DL$.
\end{abstract}

\maketitle

\section{Introduction}

In~\cite{bousfield-bool} and ~\cite{bousfield-spec}, Bousfield
introduced an equivalence relation on spectra that has turned out to
be extremely important.  Given a spectrum $E$, we define the
\emph{Bousfield class} $\bousclass{E}$ of $E$ to be the collection of
\emph{$E$-acyclic} spectra $X$, where $X$ is $E$-acyclic if and only
if $E\Smash X=0$.  Then we say that $E$ and $F$ are \emph{Bousfield
equivalent} if and only if $\bousclass{E} = \bousclass{F}$.  The
notion of Bousfield equivalence, and hence Bousfield class, plays a
major role in much of modern stable homotopy theory.

We can order the collection of Bousfield classes using reverse
inclusion.  We then have a partially ordered class associated to the
stable homotopy category, and Bousfield and others have investigated
properties of this partially ordered class.  The nilpotence theorem of
Devinatz, Hopkins, and Smith~\cite{devinatz-hopkins-smith}, for
example, is equivalent to the classification of Bousfield classes of
finite spectra~\cite{hopkins-smith}.  We recently learned that Ohkawa
has proved the surprising result that there is only a set of Bousfield
classes~\cite{ohkawa}; see
also~\cite{strickland-counting-bousfield-classes}.  He proves there
are at most $\card{2}$ Bousfield classes, where
$\card{i}=2^{\card{i-1}}$ and $\card{0}=\aleph_{0}$.  In light of this
result, the authors decided to re-examine the structure of the
partially ordered set of Bousfield classes.

The goal of this paper is to provide some kind of global understanding
of the partially ordered set $\lattice$ of Bousfield classes.  Using
Ohkawa's result, we are able to perform certain constructions in
$\lattice$, such as a greatest lower bound operation.  We also bring
to bear many methods and results from lattice theory; for instance,
the sub-partially ordered set $\DL$ of $\lattice$, which consists of
all Bousfield classes $\bousclass{X}$ for which $\bousclass{X} =
\bousclass{X \Smash X}$, is a very nice sort of distributive lattice
known as a \emph{frame}.  This has some nice consequences, and it also
leads to some interesting questions.  Much of our understanding of the
Bousfield lattice is only conjectural; we hope that the conjectures
and their implications are interesting enough to prompt further study
of this material.  There are several questions we have not addressed.
In particular, a frame such as $\DL $ has an associated topological
space.  It would be interesting to understand something about this
space, even conjecturally.  Jack Morava has asked whether this space has
a structure sheaf, probably of stable homotopy categories, associated to
it.  The stalk at $K(n)$, for example, might be the $K(n)$-local
category.  There are also many frame-theoretic properties that $\DL $
may or may not have, such as coherence.  

Here is one of the conjectures that we do discuss.  Call a Bousfield
class $\bousclass{X}$ \emph{strange} if $\bousclass{X} <
\bousclass{\HFp}$.  For instance, the Brown-Comenetz dual of the
$p$-local sphere has such a Bousfield class.  By general lattice
theory, the inclusion $\DL \hookrightarrow \lattice$ has a ``right
adjoint'' $r \mathcolon \lattice \xrightarrow{} \DL$ which is a
retraction onto $\DL$.  One can see that $r$ sends every strange
Bousfield class to $\bousclass{0}$, and also that $r$ induces a map
$r' \mathcolon \lattice / (\textup{strange}) \xrightarrow{} \DL$,
where $\lattice / (\textup{strange})$ is the quotient lattice of
$\lattice$ by the ideal of strange Bousfield classes.
Conjecture~\ref{conj-r} states that $r'$ is an isomorphism; this
implies, for example, that $\bousclass{E \Smash E} = \bousclass{E
\Smash E \Smash E}$ for all spectra $E$.  Our other main conjectures
are Conjectures~\ref{conj-A(n)-minimal}, \ref{conj-D},
\ref{conj-nothing-detects-I}--\ref{conj-everything-bigger-than-I}, and
\ref{conj-localizing-subcats}.

Here is the structure of the paper.  In Section~\ref{sec-structure},
we define Bousfield classes and the basic operations one can perform
on them: join, smash, meet, and complementation.  Next we examine
$\DL$ and its relation to $\lattice$; in particular, we note that
$\DL$ is a frame, and we construct a retraction from $\lattice$ to
$\DL$.  We also give the conjectured description of this retraction in
terms of strange Bousfield classes.  We discuss more basic structure
in Section~\ref{sec-complement}: we discuss minimal and complemented
Bousfield classes, and we recall some properties of $\BA$, the set of
complemented Bousfield classes.  For example, we recall Bousfield's
observation that $\BA$ is a Boolean algebra.  In
Section~\ref{sec-finite-acyclics} we examine spectra $X$ for which
there is a finite spectrum $F$ with $X \Smash F = 0$; we give a
conjectured classification of the Bousfield classes of such $X$.  This
provides some information about $\BA$.  In Section~\ref{sec-cBA}, we
return to the fact that $\DL$ is a frame; this allows us to construct
a complete Boolean algebra $\cBA \subseteq \DL$ which (properly)
contains $\BA$, and we give a conjectured description of $\cBA$.  Then
in the next section, we examine Bousfield classes of spectra $X$ for
which $X \Smash F \neq 0$ for all finite $F$.  This leads to a
discussion of some properties of $I$, the Brown-Comenetz dual of the
$p$-local sphere, as well as several conjectures about spectra with no
finite acyclics.  We show that these conjectures are all equivalent,
and we discuss some of their consequences.  Much of the paper to this
point suggests that the set of strange Bousfield classes, those
classes of $p$-local spectra $X$ with $\bousclass{X} <
\bousclass{\HFp}$, is interesting; in Section~\ref{sec-strange} we
examine some examples of such spectra.  We end the paper in
Section~\ref{sec-localizing} with a discussion of the partially
ordered class of localizing subcategories---recall that a subcategory
is called \emph{localizing} if it is thick and is closed under
coproducts; the main conjecture is that every localizing subcategory
is equal to the class of $E$-acyclics for some spectrum $E$.  This
conjecture has several equivalent formulations, and some deep
structural consequences.

We work $p$-locally throughout the paper, except for
Section~\ref{sec-localizing}, in which we work globally.  As in all
discussions of Bousfield classes of spectra, we work in the
stable homotopy category of spectra, as described for example
in~\cite{hovey-palmieri-strickland}.  

The authors would like to thank Dan Christensen and Neil Strickland
for many helpful discussions about Bousfield classes.

\section{Basic structure of the Bousfield lattice: $\vee$, $\Smash$,
$\meet$, and $a$}
\label{sec-structure}

In this section we discuss the basic structure of the Bousfield
lattice, including the wedge (a.k.a.~the join) $\vee$, the smash
product $\Smash$, the meet $\meet$, and the complementation operator
$a$.

We start with the definition of Bousfield equivalence and related
ideas, due to Bousfield in \cite{bousfield-bool} and
\cite{bousfield-spec}.

\begin{definition}
Let $E$, $F$, $X$, and $Z$ be spectra.
\begin{enumerate}
\item $Z$ is \emph{$E$-acyclic} if and only if $E \Smash Z = 0$.
\item The \emph{Bousfield class} of $E$, written $\bousclass{E}$, is
the collection of $E$-acyclic spectra.
\item The spectra $E$ and $F$ are \emph{Bousfield equivalent} if and
only if $\bousclass{E} = \bousclass{F}$.
\item The Bousfield classes are partially ordered by reverse
inclusion: we write $\bousclass{E} \geq \bousclass{F}$ if and only if
$E \Smash Z = 0 \Rightarrow F \Smash Z = 0$.
\item The \emph{wedge} of $\bousclass{E}$ and $\bousclass{F}$ is
defined to be $\bousclass{E} \vee \bousclass{F} = \bousclass{E \vee
F}$.  The wedge of an arbitrary set of Bousfield classes is defined
the same way.
\item Similarly, the \emph{smash product} of $\bousclass{E}$ and
$\bousclass{F}$ is defined to be $\bousclass{E} \Smash \bousclass{F} =
\bousclass{E \Smash F}$.
\item $X$ is \emph{$E$-local} if and only if $[Z,X]=0$ for all
$E$-acyclic spectra $Z$.
\end{enumerate}
\end{definition}

We denote the $p$-local sphere by $S$; then $\bousclass{S}$ is the
largest Bousfield class in this ordering, and $\bousclass{0}$ is the
smallest.  It is clear that $\bousclass{E} \vee \bousclass{F}$ is the
least upper bound, or \emph{join}, of $\bousclass{E}$ and
$\bousclass{F}$; indeed, $\bigvee \bousclass{E_{i}}$ is the join of
the set $\{\bousclass{E_{i}}\}$.

Now we recall Ohkawa's result.

\begin{theorem}[\cite{ohkawa}]
The class of Bousfield classes forms a set.
\end{theorem}

We use $\lattice$ to denote the set of Bousfield classes.

$\lattice$ is a partially ordered set in which every subset has a
least upper bound (i.e., $\lattice$ is a \emph{complete join
semilattice}).  Since there is a smallest element, then every subset
also has a greatest lower bound, or \emph{meet}, obtained by taking
the join of all the lower bounds (we are using the fact there is a
(nonempty) \emph{set} of these lower bounds, so we can in fact take
the join).  Since $\lattice$ has both finite joins and finite meets,
it is a \emph{lattice}; since it has arbitrary joins, it is a
\emph{complete} lattice.  We denote the meet of $\bousclass{X}$ and
$\bousclass{Y}$ by $\bousclass{X} \meet \bousclass{Y}$.
Unfortunately, this meet is not easily described.  In particular, we
do not know whether $\lattice$ is distributive: in other words, is
$\bousclass{X}\vee \lbp \bousclass{Y} \meet \bousclass{Z} \rbp = \lbp
\bousclass{X} \vee \bousclass{Y} \rbp \meet \lbp \bousclass{X} \vee
\bousclass{Z} \rbp$?  The meet certainly does not distribute over
infinite joins; see Example~\ref{eg-meet-infinite-join}.

In contrast, the smash product of Bousfield classes distributes over
infinite joins.  $\bousclass{X} \Smash \bousclass{Y}$ is a lower bound
for $\bousclass{X}$ and $\bousclass{Y}$, but it need not be the
greatest lower bound; for example, if $I$ is the Brown-Comenetz dual
of the sphere, then $I \Smash I = 0$ (see \cite[Lemma
2.5]{bousfield-bool} and Lemma~\ref{lem-smash-square-equal-0} below).
In general, then, we have $\bousclass{X} \Smash \bousclass{Y} \leq
\bousclass{X} \meet \bousclass{Y}$.

In any complete lattice with an operation that distributes over
infinite joins, we can define a complementation operator $a$: we
define $a\bousclass{X}$ to be the join of all $\bousclass{Y}$ such
that $\bousclass{X}\Smash \bousclass{Y}=\bousclass{0}$.  Here are some
of the basic properties of $a$, most of which are due to
Bousfield~\cite{bousfield-spec}.

\begin{lemma}\label{lem-properties-of-a}
Let $a$ be the complementation operator on the Bousfield lattice
$\lattice$.  Then $a$ has the following properties.
\begin{enumerate}
\item $\bousclass{E}\leq a\bousclass{X}$ if and only if $E\Smash X=0$.
\item $a$ is order-reversing\uc{} $\bousclass{X}\leq \bousclass{Y}$ if
and only if $a\bousclass{X}\geq a\bousclass{Y}$.
\item $a^{2}\bousclass{X}=\bousclass{X}$.
\item $\bousclass{X}\meet \bousclass{Y}=a\lbp a\bousclass{X}\vee
a\bousclass{Y} \rbp$.
\item More generally, $a$ converts arbitrary joins to meets and
arbitrary meets to joins.
\end{enumerate}
\end{lemma}

\begin{proof}
Part~(a) holds since the smash product distributes over infinite
joins, so $a\bousclass{X}\Smash \bousclass{X}=0$.  For the next part,
suppose $\bousclass{X}\leq \bousclass{Y}$.  Since
$a\bousclass{Y}\Smash \bousclass{Y}=0$, then $a\bousclass{Y}\Smash
\bousclass{X}=0$.  Hence $a\bousclass{Y}\leq a\bousclass{X}$, so $a$
is order-reversing.  The other half of part~(b) follows from part~(c).

For part~(c), it is formal to verify that $\bousclass{X}\leq
a^{2}\bousclass{X}$.  Since $a$ is order-reversing, it follows that
$a\bousclass{X}\geq a^{3}\bousclass{X}$.  Thus
$a\bousclass{X}=a^{3}\bousclass{X}$.  Now suppose $X\Smash Z=0$.  Then
$\bousclass{Z}\leq a\bousclass{X}=a\lbp a^{2}\bousclass{X} \rbp$, so
$a^{2}\bousclass{X}\Smash \bousclass{Z}=0$.  Thus $\bousclass{X}\geq
a^{2}\bousclass{X}$, completing the proof of part~(c). Parts~(d)
and~(e) are formal consequences of the other parts, given Ohkawa's
theorem.
\end{proof}

Note that not all of these properties would hold if we tried to define
$a$ using the meet instead of the smash.  

Bousfield's work predates Ohkawa's, so he had to work harder to
construct the operator $a$.  In particular, he constructs an operator
at (or closer to) the spectrum level, and shows that it descends to
give an operator on Bousfield classes.  For any spectrum $E$,
Bousfield shows in~\cite[Lemma 1.13]{bousfield-spec} that the
localizing subcategory of $E$-acyclic spectra is generated by a single
spectrum $aE$.  So for instance, a spectrum $X$ is $E$-local if and
only if $[aE,X]_{*}=0$.  The spectrum $aE$ is not well-defined, but
any other choice generates the same localizing subcategory, so in
particular has the same Bousfield class.  Thus $\bousclass{aE}$ is
well-defined.  In fact, $\bousclass{aE}=a\bousclass{E}$, since if
$Z\Smash E=0$, then $Z$ is in the localizing subcategory generated by
$aE$, so $\bousclass{Z}\leq \bousclass{aE}$.

As with the meet, it is rather difficult to compute the effect of the
operator $a$.  We will discuss it further and give some examples in
Section~\ref{sec-complement}.

Dan Christensen has pointed out that, just as one can define the meet
operation $\meet$ by $\bousclass{X} \meet \bousclass{Y} = a\lbp
a\bousclass{X} \vee a\bousclass{Y} \rbp$, one can define an operation
$\curlyvee$ by
\[
\bousclass{X} \curlyvee \bousclass{Y} = a\lbp a\bousclass{X} \Smash
a\bousclass{Y} \rbp.
\]
Then $\bousclass{X} \curlyvee \bousclass{Y} \geq \bousclass{X} \vee
\bousclass{Y}$, and this inequality may be strict; for example,
$a\bousclass{I} \curlyvee a\bousclass{I} = \bousclass{S}$, even though
$a\bousclass{I} \neq \bousclass{S}$.

\section{The retraction onto $\DL$}\label{sec-r}

The Bousfield lattice $\lattice$ is a complete lattice, but may not be
distributive; we show in Example~\ref{eg-meet-infinite-join} that the
meet does not distribute over infinite joins.  In any case, the smash
and the meet certainly do not coincide.  To get around this problem,
Bousfield introduced the sub-partially ordered set $\DL$ of $\lattice$
in~\cite{bousfield-bool}: $\DL$ consists of the Bousfield classes
$\bousclass{E}$ satisfying $\bousclass{E}=\bousclass{E} \Smash
\bousclass{E}$.  The goal of this section is to study $\DL$ and its
relationship to $\lattice$.  In particular, we point out that there is
a retraction $\lattice \xrightarrow{}\DL$, and we make some
conjectures about it.

\begin{example}
Bousfield observes in \cite{bousfield-bool} that if $E$ is a ring
spectrum or a finite spectrum, then $\bousclass{E}$ is in $\DL$.  On
the other hand, $I$, the Brown-Comenetz dual of the $p$-local sphere
is not: since $I \Smash I = 0$, then $\bousclass{I} > \bousclass{I}
\Smash \bousclass{I}$.
\end{example}

We mention the following in passing.

\begin{question}
Let $E$ be a spectrum.  Must the sequence
\[
\bousclass{E} \geq \bousclass{E}^{\Smash 2} \geq \bousclass{E}^{\Smash
3} \geq \dots
\]
stabilize?
\end{question}

See Proposition~\ref{prop-description-of-r}(e) for a conjectured
answer to this question.

\begin{lemma}\label{lem-DL-meet}
Suppose $\bousclass{E}\in \DL$, $\bousclass{E}\leq \bousclass{X}$, and
$\bousclass{E}\leq \bousclass{Y}$.  Then $\bousclass{E}\leq
\bousclass{X} \Smash \bousclass{Y}$.  
\end{lemma}

\begin{proof}
We have $\bousclass{E}=\bousclass{E}\Smash \bousclass{E}\leq
\bousclass{X}\Smash \bousclass{Y}$.  
\end{proof}

A \emph{frame} is a complete lattice in which the meet distributes
over infinite joins: $a\Smash \bigvee b_{i}=\bigvee (a\Smash b_{i})$.
For example, a topology on a space $X$ has the structure of a frame,
in which the open subsets of $X$ are ordered by inclusion.  Frames are
also called \emph{locales}, \emph{complete Heyting algebras}, or
\emph{complete Brouwerian lattices}.  They are used in categorical
topology~\cite{johnstone}, where a locale is viewed as a
generalized topological space, lattice theory~\cite{birkhoff}, and
logic~\cite{freyd-allegories}.

\begin{proposition}\label{prop-DL-lattice}
$\DL$ is a frame.  In $\DL$, the join of $\{\bousclass{X_{i}} \}$
is $\bigvee \bousclass{X_{i}}$, and the meet of $\bousclass{X}$ and
$\bousclass{Y}$ is $\bousclass{X}\Smash \bousclass{Y}$.  The inclusion
$i\mathcolon \DL \xrightarrow{}\lattice$ preserves arbitrary joins
but does not preserve meets.
\end{proposition}

\begin{proof}
Much of this is due to Bousfield \cite{bousfield-bool}, either
explicitly or implicitly.  We leave to the reader the straightforward
check that $\bigvee \bousclass{X_{i}}$ and $\bousclass{X} \Smash
\bousclass{Y}$ are in $\DL$ if all $\bousclass{X_{i}}$,
$\bousclass{X}$, and $\bousclass{Y}$ are in $\DL$.  It follows from
this that $\bigvee \bousclass{X_{i}}$ is the join of
$\{\bousclass{X_{i}} \}$, that $\DL$ is a complete lattice, and that
the inclusion $i\mathcolon \DL \xrightarrow{}\lattice$ preserves
joins.  Lemma~\ref{lem-DL-meet} implies that the meet in $\DL$ is the
smash product, and, since the smash product distributes over infinite
joins, that $\DL$ is a frame.  To see that $i$ does not preserve
meets, note that both $\bousclass{\HFp }$ and $\bigvee_{n}
\bousclass{K(n)}$ are in $\DL$.  Their smash product, and hence their
meet in $\DL$, is $0$, but their meet in $\lattice$ is at least
$\bousclass{I}$ by Proposition~\ref{prop-I-detects-finite-locals}.
\end{proof}

We can think of a complete lattice, or indeed any partially ordered
set, as a category with a unique map from $x$ to $y$ if and only if
$x\leq y$.  A complete lattice is just a partially ordered set that is
complete and cocomplete as a category; the colimit of a functor to a
lattice is the join of all of the objects in the image, and dually for
the limit.  From this point of view, an order-preserving map of
partially ordered sets corresponds to a functor on the associated
categories.  A functor between complete lattices preserves colimits if
and only if it preserves arbitrary joins.  Obviously a left adjoint
must have this property, and, for complete lattices, the converse is
true as well.  Note that, for maps of partially ordered sets $f$ and
$g$, $g$ is right adjoint to $f$ if and only if $fx\leq y$ is
equivalent to $x\leq gy$.

\begin{lemma}\label{lem-adjoints}
Suppose $f\mathcolon \mathbf{C}\xrightarrow{}\mathbf{D}$ is an
order-preserving map between complete lattices.  Then $f$ has a right
adjoint if and only if $f$ preserves arbitrary joins.  In this case, the
right adjoint of $f$ is the map $g$ defined by $gy=\lub \{x \suchthat
fx\leq y \}$.
\end{lemma}

\begin{proof}
One can easily verify that $g$ is order-preserving and $fx\leq y$ implies
$x\leq gy$.  Conversely, if $f$ preserves colimits, then $fgy=\lub \{fx
\suchthat fx\leq y \}\leq y$, so $x\leq gy$ implies $fx\leq fgy\leq y$.
\end{proof}

Johnstone proves in \cite[Theorem I.4.2]{johnstone} the (equivalent)
statement that a functor between complete lattices has a left adjoint
if and only if it preserves arbitrary meets.  
Applying Lemma~\ref{lem-adjoints}
to $\DL$, we get the following corollary, first pointed out to us by
Neil Strickland.  

\begin{corollary}\label{cor-right-adjoint}
The inclusion functor $\DL \xrightarrow{} \lattice$ has a right adjoint
$r\mathcolon \lattice \xrightarrow{}\DL$ defined by
$r\bousclass{X}=\lub \{ \bousclass{Y}\in \DL \suchthat \bousclass{Y}\leq
\bousclass{X} \}$.  The functor $r$ preserves arbitrary meets,
$r\bousclass{X}\leq \bousclass{X}$ for all $X$, and
$r\bousclass{X}=\bousclass{X}$ if $\bousclass{X}\in \DL$.
\end{corollary}

In fact $r$ preserves the smash product as well.  

\begin{lemma}\label{lem-r-preserves-smash}
The functor $r\mathcolon \lattice \xrightarrow{}\DL$ preserves the
smash product\uc{} $r\lbp \bousclass{X} \Smash
\bousclass{Y} \rbp=r\bousclass{X}\Smash
r\bousclass{Y}$.
\end{lemma}

\begin{proof}
Since $\bousclass{X}\Smash \bousclass{Y}$ is a lower bound for
$\bousclass{X}$ and $\bousclass{Y}$, $r\lbp \bousclass{X} \Smash
\bousclass{Y} \rbp$ is a lower bound for $r\bousclass{X}$ and
$r\bousclass{Y}$, so $r\lbp \bousclass{X} \Smash \bousclass{Y}\rbp
\leq r\bousclass{X}\Smash r\bousclass{Y}$.  Conversely,
$r\bousclass{X}\Smash r\bousclass{Y}\leq \bousclass{X}\Smash
\bousclass{Y}$ and $r\bousclass{X}\Smash r\bousclass{Y}\in \DL$, so
$r\bousclass{X}\Smash r\bousclass{Y}\leq r \lbp \bousclass{X} \Smash
\bousclass{Y} \rbp$.
\end{proof}

We would like to understand this map $r$ more explicitly.  We begin by
pointing out that $r$ does kill some Bousfield classes.  

\begin{lemma}\label{lem-smash-square-equal-0}
If $\bousclass{E}<\bousclass{\HFp }$ and $\bousclass{F}\leq
\bousclass{\HFp }$, then $E\Smash F=0$.  In particular, $\bousclass{E}
\Smash \bousclass{E} = 0$, so $r\bousclass{E}=0$.
\end{lemma}

\begin{proof}
We must have $E\Smash \HFp =0$, since otherwise $\bousclass{E} \geq
\bousclass{\HFp }$.  Hence $E\Smash F=0$.  In particular
$\bousclass{E}\Smash \bousclass{E}=0$.  Since $r\bousclass{E}\in \DL$,
Lemma~\ref{lem-DL-meet} implies that $r\bousclass{E}\leq
\bousclass{E}\Smash \bousclass{E}=0$, so $r\bousclass{E}=0$.
\end{proof}

By the argument in \cite[2.6]{rav-loc} (see also
Lemma~\ref{lem-properties-of-I}(c) below), $I$ does have
$\bousclass{I}<\bousclass{\HFp}$, so there are nontrivial examples of
such spectra.

\begin{definition}\label{defn-strange}
We define a spectrum $E$ to be \emph{strange} if
$\bousclass{E}<\bousclass{\HFp }$.
\end{definition}

Hence every strange spectrum is in the kernel of $r$.  We will study
some more examples of strange spectra in Section~\ref{sec-strange}.

A subset $J$ of a complete lattice $\mathbf{C}$ is called a
\emph{complete ideal} if it is closed under arbitrary joins, and if
$x\in J$ and $y\leq x$, then $y\in J$.  Every complete ideal in a
complete lattice is \emph{principal}; if we let $m$ be the join of all
the elements of $J$, then $y\in J$ if and only if $y\leq m$.  For the
complete ideal of strange spectra, we can identify the maximal element
$m$ ``explicitly.''

\begin{lemma}\label{lem-strange-is-principal}
Let $D=a\HFp \vee \HFp $.  Then the collection of
strange Bousfield classes is the principal ideal generated by
$a\bousclass{D}=a\bousclass{\HFp }\meet \bousclass{\HFp }$.
\end{lemma}

\begin{proof}
Note that $\bousclass{E}\leq a\bousclass{D}$ if and only if $E\Smash
\HFp =0$ and $E\Smash a\HFp =0$.  This second condition holds if and
only if $\bousclass{E}\leq a^{2}\bousclass{\HFp }=\bousclass{\HFp }$.
Hence $\bousclass{E}\leq a\bousclass{D}$ if and only if
$\bousclass{E}\leq \bousclass{\HFp }$ and $E\Smash \HFp =0$, that is, if
and only if $\bousclass{E}<\bousclass{\HFp }$.
\end{proof}

Given a (complete) ideal $J$ in a (complete) lattice $\mathbf{C}$, we
can define $a\equiv b \pmod{J}$ if there is some $x\in J$ such that
$a\vee x=b\vee x$.  If $J$ is principal, then $a \equiv b \pmod{J}$ if
$a \vee m = b \vee m$, where $m$ is the largest element in $J$.  The
equivalence classes under this congruence relation define a complete
lattice $\mathbf{C}/J$ (see~\cite[II.4]{birkhoff}, and note that a
complete join semilattice is a complete lattice).  The obvious
epimorphism $\mathbf{C}\xrightarrow{}\mathbf{C}/J$ preserves arbitrary
joins, and has kernel $J$.  There are often other epimorphisms with
kernel $J$; hence given a poset map $\mathbf{C} \xrightarrow{}
\mathbf{D}$ with kernel containing $J$, there may not be an induced
map $\mathbf{C}/J \xrightarrow{} \mathbf{D}$.

\begin{proposition}\label{prop-r-mods-out-strange}
Let $\strange$ be the principal ideal of strange Bousfield classes.  If
$\bousclass{X}\equiv \bousclass{Y} \pmod{\strange}$, then
$r\bousclass{X}=r\bousclass{Y}$.  
\end{proposition}

\begin{proof}
As before, we let $D=a\HFp \vee \HFp $.  Since $\strange$ is the principal
ideal generated by $a\bousclass{D}$, we have $\bousclass{X}\equiv
\bousclass{Y}\pmod{\strange}$ if and only if $\bousclass{X} \vee
a\bousclass{D}=\bousclass{Y} \vee a\bousclass{D}$.  It therefore
suffices to show that $r\lbp\bousclass{X} \vee
a\bousclass{D}\rbp=r\bousclass{X}$.  So suppose $\bousclass{Z}\in
\DL$ with $\bousclass{Z}\leq \bousclass{X}\vee a\bousclass{D}$.  Then
Lemma~\ref{lem-DL-meet} implies that $\bousclass{Z}=\lbp 
\bousclass{Z} \Smash \bousclass{X} \rbp \vee \lbp 
\bousclass{Z} \Smash a\bousclass{D} \rbp$.  Now, if $Z\Smash \HFp
=0$, then $Z\Smash aD=0$, and so $\bousclass{Z}=\bousclass{Z} \Smash
\bousclass{X}\leq \bousclass{X}$.  On the other hand, if $Z\Smash \HFp
$ is nonzero, then $(X\vee aD)\Smash \HFp $ is nonzero, so $X\Smash
\HFp $ is nonzero.  Hence $\bousclass{Z} \Smash a\bousclass{D}\leq
\bousclass{\HFp }\leq \bousclass{X}$, so $\bousclass{Z}\leq
\bousclass{X}$ in this case as well.  Thus $r\lbp
\bousclass{X} \vee a\bousclass{D} \rbp =r\bousclass{X}$ as required.
\end{proof}

It follows from Proposition~\ref{prop-r-mods-out-strange} that the
epimorphism $r\mathcolon \lattice \xrightarrow{}\DL$ factors through
an epimorphism $r'\mathcolon \lattice /\strange \xrightarrow{}\DL$,
where $\strange$ is the ideal of strange spectra.

\begin{conjecture}\label{conj-r}
The epimorphism $r'\mathcolon \lattice /\strange \xrightarrow{}\DL$ is an
isomorphism.  
\end{conjecture}

This conjecture has two parts: that $\strange$ is the kernel of $r$,
and (since epimorphisms of lattices are not determined by their
kernels) that the induced map is an isomorphism.  The conjecture has
several consequences.

\begin{proposition}\label{prop-description-of-r}
Suppose Conjecture~\ref{conj-r} holds.  Then the following properties
hold.  
\begin{enumerate}
\item $r\bousclass{E}=0$ if and only if $E$ is strange. 
\item $r$ preserves arbitrary joins. 
\item If $E\Smash \HFp \neq 0$, then $\bousclass{E}\in \DL$.  
\item $r\bousclass{E}=\bousclass{E} \Smash \bousclass{E}$.  
\item Hence $\bousclass{E}^{\Smash n} = \bousclass{E}^{\Smash (n+1)}$
when $n \geq 2$.
\end{enumerate}
\end{proposition}

\begin{proof}
The first two parts are immediate.  For part~(c), note that
Conjecture~\ref{conj-r} implies that $\bousclass{E}\equiv r\bousclass{E}
\pmod{\strange}$, so that $\bousclass{E} \vee a\bousclass{D}=r\bousclass{E}\vee
a\bousclass{D}$, where $D=a\HFp \vee \HFp $ as usual.  If $E\Smash \HFp
\neq 0$, then $\bousclass{E}\geq \bousclass{\HFp }>a \bousclass{D}$, so
$\bousclass{E} \vee a\bousclass{D}=\bousclass{E}$.  Similarly,
$r\bousclass{E}>a\bousclass{D}$, so $r\bousclass{E}\vee
a\bousclass{D}=r\bousclass{E}$.  Thus $\bousclass{E}=r\bousclass{E}$, and
so $\bousclass{E}\in \DL$.

Part~(d) is proved similarly.  We can assume that $E\Smash \HFp =0$.  Then
\begin{align*}
\bousclass{E}\Smash \bousclass{E} & = 
  \lbp \bousclass{E} \vee a\bousclass{D} \rbp \Smash 
  \lbp \bousclass{E} \vee a\bousclass{D} \rbp \\
 & = \lbp r\bousclass{E}\vee a\bousclass{D} \rbp \Smash
     \lbp r\bousclass{E}\vee a\bousclass{D} \rbp \\
 & = r\bousclass{E}\Smash r\bousclass{E} \\
 & = r\bousclass{E}.
\end{align*}
Part~(e) follows immediately.
\end{proof}

Note that, if $r$ preserves arbitrary joins, it must have a right
adjoint $r^{*}\mathcolon \DL \xrightarrow{} \lattice$.  This right
adjoint must be defined by $r^{*}\bousclass{E}=\bousclass{E}\vee
a\bousclass{D}$, where $D=a\HFp \vee \HFp $.  We can define this map
without knowing Conjecture~\ref{conj-r}, of course, but we do not know
that it preserves arbitrary meets without Conjecture~\ref{conj-r}.

Another corollary of Conjecture~\ref{conj-r} would be some understanding
of the difference between the meet and the smash in $\lattice$.  In
particular, the meet and the smash are equivalent, modulo strange
spectra.

\begin{proposition}\label{prop-meet-smash}
Suppose Conjecture~\ref{conj-r} holds.  Let $D=a\HFp \vee \HFp$, so that
$a\bousclass{D}$ is the maximum strange Bousfield class.  Then if
$\bousclass{X}$ and $\bousclass{Y}$ are arbitrary Bousfield classes, we
have
\[
\lbp \bousclass{X} \meet \bousclass{Y} \rbp \vee a\bousclass{D} =
\lbp \bousclass{X}\Smash \bousclass{Y} \rbp \vee a\bousclass{D}. 
\]
\end{proposition}

\begin{proof}
Since $r$ preserves both meets and the smash product, we have
$r\lbp \bousclass{X}\meet \bousclass{Y} \rbp=r\lbp \bousclass{X}\Smash
\bousclass{Y}\rbp$.  Conjecture~\ref{conj-r} completes the proof.  
\end{proof}

\section{More structure of $\lattice$: minimal and complemented
classes}
\label{sec-complement}

In this section we discuss minimal, maximal, and complemented
Bousfield classes.

We say that a nonzero Bousfield class $\bousclass{E}$ is
\emph{minimal} if there is no nonzero Bousfield class strictly less
than $\bousclass{E}$.  \emph{Maximal} Bousfield classes are defined
similarly.

\begin{example}
For $n \geq 0$, the \emph{$n$th Morava $K$-theory} spectrum $K(n)$
has a minimal Bousfield class---see Section~\ref{sec-finite-acyclics}.
We conjecture below (Conjecture~\ref{conj-A(n)-minimal}) that
$\bousclass{A(n)}$ is minimal when $n \geq 2$, where $A(n)$ is a
spectrum that measures the failure of the telescope conjecture; we
also conjecture (see Lemma~\ref{lem-I-is-minimal}) that
$\bousclass{I}$ is minimal, where $I$ is the Brown-Comenetz dual of
the sphere.
\end{example}

It is natural to wonder whether a given Bousfield class can be
written as the least upper bound of minimal ones, or dually, whether
a class is the greatest lower bound of maximal ones.  Since the
least upper bound has a much more convenient description, we will
focus on minimal Bousfield classes.  Using the complementation
operator $a$, one can easily check that $\bousclass{X}$ is minimal if
and only if $a\bousclass{X}$ is maximal.  

Although we have referred to $a$ as the complementation operator, it is not
always the case that $a\bousclass{X}\vee \bousclass{X}=\bousclass{S}$;
when this happens, we say that $\bousclass{X}$ is \emph{complemented}.
One can easily check that if there is a Bousfield class
$\bousclass{Y}$ so that $\bousclass{X}\vee
\bousclass{Y}=\bousclass{S}$ and $\bousclass{X}\Smash
\bousclass{Y}=\bousclass{0}$, then $\bousclass{Y}=a\bousclass{X}$.
This is the reason for the term ``complemented.''  We also define
$\bousclass{X}$ to be \emph{$\meet$-complemented} if there is a
Bousfield class $\bousclass{Y}$ so that $\bousclass{X} \meet
\bousclass{Y} = \bousclass{0}$ and $\bousclass{X} \vee \bousclass{Y}
= \bousclass{S}$.

Now we note that we should only have made one definition.

\begin{proposition}\label{prop-meet-complement}
$\bousclass{X}$ is $\meet$-complemented if and only if $\bousclass{X}$
is complemented.  If these conditions hold, then the
$\meet$-complement of $\bousclass{X}$ is $a\bousclass{X}$.  
\end{proposition}

\begin{proof}
Since $\bousclass{X} \Smash \bousclass{Y} \leq \bousclass{X} \meet
\bousclass{Y}$, we see that if $\bousclass{X}$ is
$\meet$-complemented, then $\bousclass{X}$ is complemented, with the
same complement.  Conversely, suppose that $a\bousclass{X}\vee
\bousclass{X}=\bousclass{S}$.  Then
\[
a\bousclass{X}\meet \bousclass{X}=a\lbp a^{2}\bousclass{X}\vee
a\bousclass{X}\rbp = a\bousclass{S}=\bousclass{0},
\]
so $a\bousclass{X}$ is the $\meet $-complement of $\bousclass{X}$. 
\end{proof}

The collection of all complemented Bousfield classes is denoted $\BA$.
Here are some of the basic properties of $\BA$; these are all
due to Bousfield~\cite{bousfield-bool}.  

\begin{lemma}\label{lem-properties-of-BA}
Suppose that $\bousclass{X}$ and $\bousclass{Y}$ are in $\BA$, and
$\bousclass{E}$ is an arbitrary Bousfield class.  Then\uc 
\begin{enumerate}
\item $\bousclass{E}=\lbp \bousclass{E} \Smash \bousclass{X} \rbp \vee
\lbp \bousclass{E} \Smash a\bousclass{X} \rbp$.
\item $\bousclass{E}\leq \bousclass{X}$ if and only if
$\bousclass{E}=\bousclass{E} \Smash \bousclass{X}$.  
\item $\bousclass{X}\meet \bousclass{Y}=\bousclass{X}\Smash \bousclass{Y}$.  
\item Hence $\BA \subseteq \DL$.
\item $\bousclass{X} \Smash \bousclass{Y}$ is in $\BA$, and $a \lbp
\bousclass{X} \Smash \bousclass{Y} \rbp =a\bousclass{X}\vee a\bousclass{Y}$.  
\item $\bousclass{X} \vee \bousclass{Y}$ is in $\BA$, and $a\lbp
\bousclass{X} \vee \bousclass{Y} \rbp =a\bousclass{X}\Smash a\bousclass{Y}$.  
\item $\BA$ is a Boolean algebra.  
\end{enumerate}
\end{lemma}

(Recall that a \emph{Boolean algebra} is a distributive lattice in
which every element has a complement.)

\begin{proof}
For the first part, use the identity $\bousclass{E}=\bousclass{E}
\Smash \bousclass{S}=\bousclass{E}\Smash \lbp \bousclass{X}\vee
a\bousclass{X}\rbp$.  The second part then follows immediately.  For
part~(c), suppose $\bousclass{E}\leq \bousclass{X}$ and
$\bousclass{E}\leq \bousclass{Y}$.  Then
$\bousclass{E}=\bousclass{E}\Smash \bousclass{X}=\lbp \bousclass{E}
\Smash \bousclass{Y} \rbp \Smash \bousclass{X}$.  Hence
$\bousclass{E}\leq \bousclass{X} \Smash \bousclass{Y}$.  Part~(d) is
clear.  For part~(e), note that $a\lbp \bousclass{X}\Smash
\bousclass{Y}\rbp =a\lbp \bousclass{X}\meet \bousclass{Y}\rbp
=a\bousclass{X}\vee a\bousclass{Y}$.  Furthermore,
\begin{align*}
\bousclass{S} & = a\bousclass{X}\vee \bousclass{X} \\
 & = a\bousclass{X}\vee \lbp \bousclass{X} \Smash \bousclass{Y} \rbp
\vee \lbp \bousclass{X}\Smash a\bousclass{Y}\rbp \\
 & \leq \lbp \bousclass{X}\Smash \bousclass{Y}\rbp \vee
a\bousclass{X}\vee a\bousclass{Y}.
\end{align*}
Thus $\bousclass{X} \Smash \bousclass{Y}$ is complemented, as
required.  The proof of part~(f) is similar, and part~(g) follows
immediately from the preceding parts.
\end{proof}

\begin{example}
Bousfield shows in \cite{bousfield-bool} that if $F$ is a finite
spectrum, then $\bousclass{F}$ is in $\BA$.  He also notes that
$\bousclass{\HZ}$ is not in $\BA$; in particular, the inclusion $\BA
\subset \DL$ is proper.  We show in Section~\ref{sec-finite-acyclics}
that $\bousclass{K(n)}$ and $\bousclass{A(n)}$ are in $\BA$.
\end{example}

The structure theory of infinite Boolean algebras is considerably more
complicated than the structure theory of finite Boolean algebras.  In
particular, $\BA$ is not closed under infinite joins (see
Corollary~\ref{cor-complemented}), and so is certainly not isomorphic to
the complete Boolean algebra of all subsets of some infinite set.  The
simplest infinite Boolean algebra that is not complete is the Boolean
algebra of all finite and cofinite subsets of an infinite set.

We have noted that every finite spectrum is complemented; some other
examples of complemented spectra are provided by smashing
localizations.  Recall that every spectrum $E$ determines a Bousfield
localization functor $L_{E}$, as described in \cite{bousfield-spec}.
If $E$ and $F$ are Bousfield equivalent, then the functors $L_{E}$ and
$L_{F}$ are equal---Bousfield localization only depends on the
Bousfield class of the spectrum.  We say that a Bousfield class
$\bousclass{E}$ is \emph{smashing} if the natural map $L_{E}S \Smash X
\xrightarrow{} L_{E} X$ is an equivalence.  Ravenel proves the
following in \cite[1.31]{rav-loc}.

\begin{proposition}\label{prop-smashing-complemented}
Every smashing Bousfield class $\bousclass{E}$ is complemented, with
complement given by the fiber $A_{E}S$ of $S \xrightarrow{} L_{E} S$.
\end{proposition}

\begin{proof}
For a general Bousfield localization functor $L_{E}$, we have
$\bousclass{S}=\bousclass{L_{E}S}\vee \bousclass{A_{E}S}$.  Because
$L_{E}$ is smashing, we have $L_{E}S\Smash A_{E}S=L_{E}A_{E}S=0$.
\end{proof}

\section{Bousfield classes with finite acyclics}\label{sec-finite-acyclics}

In this section we give a brief summary of what is known about
Bousfield classes which contain finite spectra; this leads to
information about the Boolean algebra $\BA$.  Details can be found
in~\cite{hovey-chrom-split}.

As above, we denote the
($p$-local) sphere by $S$; we write $M(p)$ for the mod $p$ Moore spectrum.
A generic finite spectrum of type $n$ will be denoted by $F(n)$; then
any choice for $F(n)$ generates the same thick subcategory
$\cat{C}_{n}$, by the thick subcategory theorem of
Hopkins-Smith~\cite{hopkins-smith,rav-nilp}.  In particular, the
Bousfield class of $F(n)$ is well-defined.  Any $F(n)$ has an
essentially unique $v_{n}$-self map whose cofiber is an $F(n+1)$ and
whose telescope we will denote by $T(n)$.  The Bousfield class of
$T(n)$ is also well-defined.

By repeated use of \cite[1.34]{rav-loc}, we have a Bousfield class
decomposition
\[
\bousclass{S} = \bousclass{T(0)} \vee \bousclass{T(1)} \vee \dots \vee
\bousclass{T(n-1)} \vee \bousclass{F(n)}.
\]
Furthermore, $T(i)\Smash T(j)=0$ unless $i=j$, and $T(i)\Smash
F(n)=0$ for $i<n$.  

It follows that localization with respect to $T(0)\vee T(1)\vee \dots
\vee T(n-1)$, written $L_{n-1}^{f}$, is smashing and that its kernel
is precisely the localizing subcategory generated by
$F(n)$---see~\cite{miller-finite}.  By the above decomposition (see
also Proposition~\ref{prop-smashing-complemented}), $\bousclass{F(n)}$
is complemented with complement $\bousclass{L_{n-1}^{f}S}$; in other
words, we have $\bousclass{S}=\bousclass{L_{n-1}^{f}S}\vee
\bousclass{F(n)}$, and $F(n)\Smash L_{n-1}^{f}S=0$.

Given a spectrum $E$, we say that $E$ \emph{has a finite acyclic} if
there is a nontrivial finite spectrum $X$ such that $E\Smash X=0$.  In
this case, the thick subcategory theorem says that the collection of
finite $E$-acyclics is $\cat{C}_{n}$ for some finite $n$, and we
have $\bousclass{E}\leq \bousclass{L_{n-1}^{f}S}$.

The Morava $K$-theory spectra $K(n)$ play an important role here.
They are known to be \emph{field spectra}, so that $K(n) \Smash E$ is
a wedge of suspensions of $K(n)$ for any $E$.  The telescope
conjecture, recently proved to be false for $n=2$ by Ravenel,
asserts that $\bousclass{T(n)}=\bousclass{K(n)}$.  If this were true,
then for any $E$ with a finite acyclic, we would have
\[
\bousclass{E}=\bigvee_{n} \bousclass{E \Smash K(n)} 
  =\bigvee _{\{n\suchthat E\Smash K(n)\neq 0\}} \bousclass{K(n)}.
\]

The failure of the telescope conjecture is measured by the fiber
$A(n)$ of the natural map $T(n)\xrightarrow{}L_{K(n)}T(n)$.  Once
again, $A(n)$ is well-defined up to Bousfield class.  With a little
work, we have $\bousclass{A(n)} \vee \bousclass{K(n)} =
\bousclass{T(n)}$; clearly $A(n)\Smash K(n)=0$.  It follows easily
from this that $\bousclass{K(n)}$ and $\bousclass{A(n)}$ are both
complemented, as of course is $\bousclass{T(n)}$.  Since $K(n)$ is a
complemented field spectrum, then $\bousclass{K(n)}$ is minimal, by
\cite[3.7.3]{hovey-palmieri-strickland}.

The spectrum $A(n)$ is rather odd, as for example $\bousclass{A(n)}
\Smash \bousclass{A(n)}=\bousclass{A(n)}$, yet $BP\Smash A(n)=0$.  So,
for instance, $A(n)$ is not (Bousfield equivalent to) a nonzero ring
spectrum.  As far as detecting finite spectra goes, $A(n)$ behaves as
$K(n)$ and $T(n)$ do:
\[
\bousclass{A(n)} \Smash \bousclass{F(i)} = \begin{cases}
\bousclass{A(n)} & \text{if $i \leq n$}, \\
0 & \text{if $i>n$}.
\end{cases}
\]
Other than this, very little is known about $A(n)$.  Since the
telescope conjecture fails when $n=2$, it seems likely that it fails
for all $n \geq 2$, in which case $A(n)$ is nonzero when $n \geq 2$.
We make the following conjectures.  The first is a replacement, of
sorts, for the telescope conjecture; it says that, although the
telescope conjecture is false, the spectra $A(n)$ that measure its
failure behave as well as possible.

\begin{conjecture}\label{conj-A(n)-minimal}
If $n\geq 2$, $\bousclass{A(n)}$ is a minimal nonzero Bousfield class.
Hence, if $E$ has a finite acyclic, then $E$ is Bousfield equivalent
to a finite wedge of spectra $K(n)$ and $A(n)$\textup{;} in
particular,
\[
\bousclass{E}=\bigvee_{\{n\suchthat E\Smash K(n)\neq 0\}}
\bousclass{K(n)} \vee \bigvee_{\{n\suchthat E\Smash A(n)\neq 0\}}
\bousclass{A(n)}.
\]
\end{conjecture}

Note that each of the wedges here is finite.
This would mean that there are only countably many Bousfield classes
with a finite acyclic.  We also have the following proposition, whose
proof is immediate.  

\begin{proposition}\label{prop-finite-acyclics-are-complemented}
Suppose Conjecture~\ref{conj-A(n)-minimal} holds.  Then every Bousfield
class with a finite acyclic is complemented. 
\end{proposition}

\section{The complete Boolean algebra of spectra}\label{sec-cBA}

We have seen that the sublattice $\DL$ of the Bousfield lattice is a
frame, and that the retraction map $r\mathcolon \lattice
\xrightarrow{}\DL$ preserves arbitrary meets.  We have conjectured
that $r$ preserves arbitrary joins.  We have not discussed how $r$
behaves with respect to complements, however, and we do so in this
section.  We also explore the relationship between $\DL$ and its
sub-poset $\BA$.

\begin{definition}\label{defn-complement-DL}
Define the \emph{complement} operation $A\mathcolon \DL
\xrightarrow{}\DL$ by
\begin{align*}
\DL & \xrightarrow{\ A\ } \DL, \\
\bousclass{X} & \mapsto r(a \bousclass{X}).
\end{align*}
\end{definition}

Then we have the following straightforward lemma, whose proof we leave
to the reader.  

\begin{lemma}\label{lem-prop-of-A}
\begin{enumerate}
\item If $\bousclass{X}$ and $\bousclass{Y}$ are in $\DL$, then
$\bousclass{Y}\leq A\bousclass{X}$ if and only if $Y\Smash X=0$.  
In other words, $A\bousclass{X} = \bigvee \{\bousclass{Y} \in \DL
\suchthat \bousclass{Y} \Smash \bousclass{X} = 0\}$.
\item $A$ is order-reversing\uc{} if $\bousclass{X}\leq \bousclass{Y}$
in $\DL$, then $A\bousclass{X}\geq A\bousclass{Y}$.
\item If $\bousclass{X}\in \DL$, then $\bousclass{X}\leq
A^{2}\bousclass{X}$ and $A\bousclass{X}=A^{3}\bousclass{X}$.  
\item $A$ converts arbitrary joins to meets\uc{} if
$\bousclass{X_{i}}$ is in $\DL$ for all $i$, then $A\lbp \bigvee
\bousclass{X_{i}} \rbp$ is the meet of the $A\bousclass{X_{i}}$.
\end{enumerate}
\end{lemma}

Note that this lemma actually holds in any frame, and the complement
operator is well-known in the theory.  See~\cite[V.11]{birkhoff}, for
example.  We will recall some of this theory in the results below for
the reader's convenience.

Also note that $A$ does not convert meets to joins.  For example, let
$X=\bigvee _{n}K(n)$ and let $Y=\HFp $.  Then $X$ and $Y$ are both in
$\DL$, and $X\Smash Y=0$, and thus $A\lbp \bousclass{X} \Smash
\bousclass{Y} \rbp=\bousclass{S}$.  On the other hand, by the
computations in Example~\ref{eg-meet-infinite-join}, we have
\[
A\bousclass{X}\vee A\bousclass{Y}\leq a\bousclass{X}\vee a\bousclass{Y}=
a\lbp \bousclass{X} \meet \bousclass{Y}\rbp \leq a\bousclass{I}<\bousclass{S}.
\]
Of course, we do have $A\lbp \bousclass{X} \Smash \bousclass{Y} \rbp
\geq A\bousclass{X}\vee A\bousclass{Y}$ for any $\bousclass{X}$ and
$\bousclass{Y}$ in $\DL$.

This argument also implies that $A^{2}$ is not the identity---indeed, if
$A^{2}$ were the identity, one can check that $A$ would have to convert
meets to joins.  However, we do not know a specific spectrum $X$ in $\DL
$ for which $A^{2}\bousclass{X}\neq \bousclass{X}$.  Given
Conjecture~\ref{conj-r}, $a\bousclass{I}$ is in $\DL$ by
Proposition~\ref{prop-description-of-r}(c), and 
$A\lbp a\bousclass{I}\rbp=r\bousclass{I}=0$, so
$A^{2}\lbp a\bousclass{I}\rbp =\bousclass{S}$.

\begin{definition}\label{defn-closed}
A Bousfield class $\bousclass{X}$ is \emph{closed} if
$\bousclass{X}\in \DL$ and $A^{2}\bousclass{X}=\bousclass{X}$.  The
sub-partially ordered set of $\DL$ consisting of the closed elements
is denoted $\cBA $.
\end{definition}

Note that every Bousfield class of the form $A\bousclass{X}$ is closed,
by Lemma~\ref{lem-prop-of-A}(c).  We have the following
theorem, which again holds in considerably more generality than we
state it; see \cite[V.10--11]{birkhoff} for the general approach.

\begin{theorem}\label{thm-cBA}
The sub-poset $\cBA $ of $\DL$ is closed under arbitrary meets, and
therefore is a complete lattice.  The join in $\cBA $ of
$\{\bousclass{X_{i}} \}$ is $A^{2}\lbp \bigvee \bousclass{X_{i}}\rbp$.
Every element in $\cBA $ is complemented, so $\cBA $ is in fact a
complete Boolean algebra.  The inclusion $\cBA \xrightarrow{}\DL$
preserves arbitrary meets, and its left adjoint is given by
$A^{2}\mathcolon \DL \xrightarrow{}\cBA $.
\end{theorem}

We will write the join in $\cBA$ as $\vee_{\cBA}$.

\begin{proof}
Note that $A^{2}$ is order-preserving.  Thus, if we denote by $\bigwedge
_{i}\bousclass{X_{i}}$ the meet in $\DL$ of $\{\bousclass{X_{i}} \}$, we
have $\bigwedge _{i}\bousclass{X_{i}} \leq A^{2} \lbp
\bigwedge_{i}\bousclass{X_{i}}\rbp \leq \bigwedge
_{i}A^{2}\bousclass{X_{i}}$.  In particular, if each $\bousclass{X_{i}}$
is closed, so is $\bigwedge _{i}\bousclass{X_{i}}$.  So $\cBA $ is
closed under arbitrary meets, and hence is a complete lattice, with the
join defined to be the meet of all upper bounds.

Now, certainly $A^{2}\lbp \bigvee_{i} \bousclass{X_{i}} \rbp$ is closed
and is an upper bound for $\{\bousclass{X_{i}} \}$.  If $\bousclass{Z}$
is closed and an upper bound for $\{\bousclass{X_{i}}\}$, we have
$\bousclass{Z}=A^{2}\bousclass{Z}\geq A^{2}\bousclass{\bigvee
_{i}X_{i}}$, so the join in $\cBA $ is as claimed.  One can easily check
that $A^{2}$ is the left adjoint to the inclusion.

It remains to show that an arbitrary element $\bousclass{X}$ of $\cBA $
is complemented in $\cBA $.  To see this, note that $\bousclass{X}\Smash
A\bousclass{X}=0$, and
\[
A^{2}\lbp \bousclass{X}\vee A\bousclass{X} \rbp=
A \lbp A\bousclass{X}\Smash A^{2}\bousclass{X}\rbp=A(0)=\bousclass{S},
\]
since
$A$ converts joins to meets.  Thus $A\bousclass{X}$ is the complement of
$\bousclass{X}$ in $\cBA $, so $\cBA $ is a complete Boolean algebra.
\end{proof}

This theorem explains our choice of symbol $\cBA $.  Note that a
complete Boolean algebra need not be isomorphic to the lattice of
subsets of a set.

Note that, if $\bousclass{X}$ is already complemented in the Bousfield
lattice, so that $\bousclass{X}\in \BA$, then certainly
$A^{2}\bousclass{X}=\bousclass{X}$, so $\BA$ is a subBoolean algebra
of $\cBA$.  Of course, the inclusion $\BA \subset \cBA$ is proper,
because $\cBA$ is complete and $\BA$ is not.  Also, the lattice $\cBA$
is not a sublattice of the Bousfield lattice: the meets and joins are
different in the two sets.

We now investigate how $A$ and $A^{2}$ behave on meets.  The following
lemma appears in~\cite[V.11]{birkhoff}; we reproduce its proof for the
reader's convenience.  

\begin{lemma}\label{lem-A-meets}
Suppose $\bousclass{X}$ and $\bousclass{Y}$ are in $\DL$.  Then 
\begin{enumerate}
\item $A\lbp \bousclass{X}\Smash \bousclass{Y}\rbp =A\lbp
A^{2}\bousclass{X}\Smash A^{2}\bousclass{Y}\rbp$.
\item $A$ converts meets to joins in $\cBA $\uc{} that is, $A\lbp
\bousclass{X}\Smash \bousclass{Y}\rbp =A^{2}\lbp A\bousclass{X}\vee
A\bousclass{Y}\rbp$.
\item $A^{2}$ preserves finite meets\uc{} that is, $A^{2}\lbp
\bousclass{X}\Smash \bousclass{Y}\rbp=A^{2}\bousclass{X}\Smash
A^{2}\bousclass{Y}$.
\end{enumerate}
\end{lemma}

\begin{proof}
Certainly $A\lbp \bousclass{X}\Smash \bousclass{Y} \rbp \geq A \lbp
A^{2}\bousclass{X}\Smash A^{2}\bousclass{Y}\rbp$.  Conversely, suppose
$\bousclass{Z}\leq A \lbp \bousclass{X}\Smash \bousclass{Y} \rbp$, so that
$Z\Smash X\Smash Y=0$.  It suffices to show that
\[
\bousclass{Z'}= \bousclass{Z}\Smash
A^{2}\bousclass{X}\Smash A^{2}\bousclass{Y}=0
\]
as well.  To see this, note that $\bousclass{Z'}\Smash
\bousclass{X}\Smash \bousclass{Y}=0$, so $\bousclass{Z'}\Smash
\bousclass{X}\leq A\bousclass{Y}$.  On the other hand,
$\bousclass{Z'}\Smash \bousclass{X}\leq A^{2}\bousclass{Y}$ by
definition.  Thus $\bousclass{Z'}\Smash \bousclass{X}\leq
A\bousclass{Y}\Smash A^{2}\bousclass{Y}=0$.  Hence $\bousclass{Z'}\leq
A\bousclass{X}$.  Since $\bousclass{Z'}\leq A^{2}\bousclass{X}$ by
definition, we have $\bousclass{Z'}\leq A\bousclass{X}\Smash
A^{2}\bousclass{X}=0$.

Part~(b) follows from part~(a), since $A$ converts joins to meets, so that
\[
A^{2}\lbp A\bousclass{X}\vee \bousclass{Y}\rbp =
  A\lbp A^{2}\bousclass{X}\Smash A^{2}\bousclass{Y}\rbp.
\]
Similarly, part~(c) follows from part~(b), since 
\[
A^{2}\lbp \bousclass{X}\Smash \bousclass{Y}\rbp =A^{3}\lbp A\bousclass{X}\vee
A\bousclass{Y}\rbp = A\lbp A\bousclass{X}\vee A\bousclass{Y}\rbp =
A^{2}\bousclass{X}\Smash A^{2}\bousclass{Y}.  \qed
\]
\renewcommand{\qed}{}\end{proof}

This lemma allows us to understand the map $A^{2}\mathcolon \DL
\xrightarrow{}\cBA $.  

\begin{definition}\label{defn-dense}
A Bousfield class $\bousclass{Z}$ is said to be \emph{dense} if
$\bousclass{Z}\in \DL$ and $A^{2}\bousclass{Z}=\bousclass{S}$.  
\end{definition}

The following theorem is a special case of Theorem~V.26
of~\cite{birkhoff}, where it is attributed to Glivenko.

\begin{theorem}\label{thm-glivenko}
For $\bousclass{X}$ and $\bousclass{Y}$ in $\DL$,
$A^{2}\bousclass{X}=A^{2}\bousclass{Y}$ if and only if there is a
dense Bousfield class $\bousclass{Z}$ such that $\bousclass{X} \Smash
\bousclass{Z}=\bousclass{Y} \Smash \bousclass{Z}$.
\end{theorem}

\begin{proof}
First suppose there is a dense $\bousclass{Z}$ such that $\bousclass{X}
\Smash \bousclass{Z}=\bousclass{Y} \Smash \bousclass{Z}$.  Then $A^{2}
\lbp \bousclass{X} \Smash \bousclass{Z} \rbp =A^{2}\lbp \bousclass{Y}
\Smash \bousclass{Z} \rbp$.  But since $A^{2}$ preserves finite meets,
this means that $A^{2}\bousclass{X}\Smash
A^{2}\bousclass{Z}=A^{2}\bousclass{Y}\Smash A^{2}\bousclass{Z}$.  Since
$A^{2}\bousclass{Z}=\bousclass{S}$, this means
$A^{2}\bousclass{X}=A^{2}\bousclass{Y}$.

Conversely, suppose $A^{2}\bousclass{X}=A^{2}\bousclass{Y}$.  Let
$\bousclass{Z}=\lbp \bousclass{X}\vee A\bousclass{Y} \rbp \Smash \lbp
A\bousclass{X}\vee \bousclass{Y}\rbp $.  Then one can easily check that
$\bousclass{X}\Smash \bousclass{Z}=\bousclass{Y}\Smash \bousclass{Z}$,
so it remains to prove that $\bousclass{Z}$ is dense.  To see this, note
that $A^{2} \mathcolon \DL \xrightarrow{} \cBA$ preserves joins, so
\[
A^{2}\lbp \bousclass{X}\vee A\bousclass{Y} \rbp
 = A^{2} \lbp A^{2}\bousclass{X}\vee_{\cBA} A^{3}\bousclass{Y} \rbp 
 = A^{2} \lbp A^{2}\bousclass{Y}\vee_{\cBA} A\bousclass{Y}\rbp 
 =\bousclass{S}.
\]
as required. 
\end{proof}



Theorem~\ref{thm-glivenko} leads us to consider the dense Bousfield
classes.

\begin{lemma}\label{lem-dense}
Let $\bousclass{D}=a\bousclass{\HFp} \vee \bousclass{\HFp}$.  If $Z$
is in $\DL$ and $\bousclass{Z}\geq \bousclass{D}$, then
$\bousclass{Z}$ is dense.  Conversely, if Conjecture~\ref{conj-r}
holds, then an arbitrary Bousfield class $\bousclass{Z} \in \lattice$
is dense if and only if $\bousclass{Z}\geq \bousclass{D}$.
\end{lemma}

\begin{proof}
If $\bousclass{Z}\geq \bousclass{D}$, then
$A\bousclass{Z}=ra\bousclass{Z}\leq ra\bousclass{D}=0$, since
$a\bousclass{D}$ is the maximum strange Bousfield class.  Hence
$A^{2}\bousclass{Z}=\bousclass{S}$, so $\bousclass{Z}$ is dense.  If
Conjecture~\ref{conj-r} holds, then any $Z$ with $\bousclass{Z}\geq
\bousclass{D}$ is automatically in $\DL$, so we can drop
that hypothesis.  Furthermore, if $\bousclass{Z}$ is dense, then
$A\bousclass{Z}=A^{3}\bousclass{Z}=A\bousclass{S}=0$, so
$ra\bousclass{Z}=0$.  Given Conjecture~\ref{conj-r}, we can conclude
that $a\bousclass{Z}$ is strange, and so $a\bousclass{Z}\leq
a\bousclass{D}$.  Thus $\bousclass{Z}\geq \bousclass{D}$.
\end{proof}

The following corollary is an immediate consequence of
Lemma~\ref{lem-dense} and Theorem~\ref{thm-glivenko}.  

\begin{corollary}\label{cor-char-of-A2}
Let $\bousclass{D}=a\bousclass{\HFp} \vee \bousclass{\HFp }$.  Suppose
Conjecture~\ref{conj-r} holds.  Then for any $\bousclass{X}$ and
$\bousclass{Y}$ in $\DL$, 
$A^{2}\bousclass{X}=A^{2}\bousclass{Y}$ if and only if
$\bousclass{X}\Smash \bousclass{D}=\bousclass{Y}\Smash \bousclass{D}$.
\end{corollary}

This corollary suggests that a characterization of $\cBA $ can be
obtained from $\bousclass{D}=a\bousclass{\HFp} \vee \bousclass{\HFp}$.
Let $\locallattice $ denote the sub-partially ordered set of
$\lattice$ consisting of all elements of the form $\bousclass{X}\Smash
\bousclass{D}$.  Then $\locallattice $ is closed under arbitrary
joins, and so is a complete lattice.  The inclusion $\locallattice
\xrightarrow{}\lattice$ preserves those arbitrary joins, so has a
right adjoint $\lattice \xrightarrow{}\locallattice$; this right
adjoint takes $\bousclass{X}$ to 
\[
\bigvee \{\bousclass{Z} \in \locallattice \suchthat
\bousclass{Z}\leq \bousclass{X} \}.
\]
If Conjecture~\ref{conj-r} holds, then $D\in \DL $, so
$\bousclass{Y}\Smash \bousclass{D}\leq \bousclass{X}$ implies that
$\bousclass{Y}\Smash \bousclass{D}\leq \bousclass{X}\Smash
\bousclass{D}$.  Thus, assuming Conjecture~\ref{conj-r}, the right
adjoint $\lattice \xrightarrow{}\locallattice $ is just given by
smashing with $\bousclass{D}$.  Smashing with $D$ preserves arbitrary
joins, so has a right adjoint $\locallattice \xrightarrow{}\lattice$
as well.  This right adjoint takes $\bousclass{X} \in \locallattice$
to the largest $\bousclass{Y}$ such that $\bousclass{Y}\Smash
\bousclass{D}=\bousclass{X}$.

\begin{lemma}\label{lem-D-in-DL}
Suppose Conjecture~\ref{conj-r} holds.  Then $\locallattice \subseteq
\DL$.  
\end{lemma}

\begin{proof}
By Conjecture~\ref{conj-r}, we have $\bousclass{X} \vee
a\bousclass{D}=r\bousclass{X}\vee a\bousclass{D}$.  Thus 
\[
\bousclass{D}\Smash \bousclass{X} 
= \bousclass{D}\Smash \lbp \bousclass{X}\vee a\bousclass{D}\rbp 
= \bousclass{D}\Smash \lbp r\bousclass{X}\vee a\bousclass{D}\rbp 
= \bousclass{D}\Smash r\bousclass{X}.
\]
Furthermore, we have $r\lbp \bousclass{D} \Smash \bousclass{X}\rbp
=r\bousclass{D}\Smash r\bousclass{X} = \bousclass{D}\Smash
r\bousclass{X}$, since Conjecture~\ref{conj-r} also implies that $D$
is in $\DL$.  Thus $r\lbp \bousclass{D} \Smash \bousclass{X}\rbp
=\bousclass{D}\Smash \bousclass{X}$, so $\bousclass{D}\Smash
\bousclass{X}$ is in $\DL$ for all $\bousclass{X}$.
\end{proof}

\begin{theorem}\label{thm-cBA-aH-and-H}
Suppose Conjecture~\ref{conj-r} holds.  Then $A^{2}\mathcolon \DL
\xrightarrow{}\cBA $ factors through the epimorphism $\bousclass{D}\Smash
(-)\mathcolon \DL \xrightarrow{}\locallattice $ to define an
isomorphism $F\mathcolon \locallattice \xrightarrow{}\cBA $.  
\end{theorem}

\begin{proof}
We define $F\lbp \bousclass{D} \Smash \bousclass{X} \rbp
=A^{2}\bousclass{X}$.  By Corollary~\ref{cor-char-of-A2}, $F$ is
well-defined, injective, and order-preserving.  On the other hand, $F$
is obviously surjective since $A^{2}$ is.
\end{proof}

Naturally, we would like a better description of $\locallattice $,
in light of Theorem~\ref{thm-cBA-aH-and-H}.  See
Conjecture~\ref{conj-A(n)-minimal} for a related result.

\begin{conjecture}\label{conj-D}
We have
\[
\bousclass{D}
 = \bigvee _{n\geq 0}\bousclass{K(n)} \vee
   \bigvee _{n\geq 2}\bousclass{A(n)}\vee \bousclass{\HFp }.
\]
\end{conjecture}

Note that $\bousclass{K(n)}\leq a\bousclass{\HFp }$ and
$\bousclass{A(n)}\leq a\bousclass{\HFp }$ for all $n$, so the $\geq$
half of the equality in Conjecture~\ref{conj-D} holds.

By the definition of $D$ and the computations in
Section~\ref{sec-finite-acyclics}, the conjecture is equivalent to
the following:
\[
\bousclass{D} =  a\bousclass{\HFp} \vee \bousclass{\HFp}
 = \bigvee_{n\geq 0} \bousclass{T(n)} \vee \bousclass{\HFp}.
\]

The following proposition completes our conjectural identification of
$\cBA $ up to isomorphism.

\begin{proposition}\label{prop-cBA-description}
Suppose Conjectures~\ref{conj-r}, \ref{conj-A(n)-minimal} and
\ref{conj-D} hold.  Then $\cBA $ is isomorphic to the complete Boolean
algebra generated by the atoms $\bousclass{K(n)}$ for $n\geq 0$,
$\bousclass{A(n)}$ for $n\geq 2$, and $\bousclass{\HFp }$.
\end{proposition}

This isomorphism is given by applying $A^{2}$, so to actually identify
$\cBA $ we need to understand the behavior of $A^{2}$.  

\begin{proposition}\label{prop-H-is-not-closed}
Suppose Conjecture~\ref{conj-r}, \ref{conj-A(n)-minimal} and
\ref{conj-D} hold.  Then every subwedge of $\bigvee _{n\geq
0}\bousclass{K(n)}\vee \bigvee _{n\geq 2}\bousclass{A(n)}$ is closed.
However, $A^{2}\bousclass{\HFp }\neq \bousclass{\HFp }$.  
\end{proposition}

\begin{proof}
Let $\bousclass{E}$ denote an arbitrary subwedge of $\bousclass{D}$ such
that $\bousclass{E}\Smash \bousclass{\HFp }=\bousclass{0}$.  Let
$\bousclass{E'}$ denote the complementary subwedge of $\bousclass{D}$.
We will show that $A\bousclass{E'}=\bousclass{E}$, so that
$\bousclass{E}$ is closed.  It is clear that $\bousclass{E}\leq
A\bousclass{E'}$, since $\bousclass{E}\Smash
\bousclass{E'}=\bousclass{0}$ and $\bousclass{E}\in \DL $.  On the other
hand, $\bousclass{E'}\geq \bousclass{\HFp }$, so $A\bousclass{E'}\leq
A\bousclass{\HFp }\leq \bousclass{D}$.  Since $A\bousclass{E'}\in \DL $,
it follows that 
\[
A\bousclass{E'} =A\bousclass{E'}\Smash \bousclass{D}
\]
and so $A\bousclass{E'}$ is a subwedge of $\bousclass{D}$.  This
subwedge cannot contain any term in $\bousclass{E'}$, so we must have
$A\bousclass{E'}=\bousclass{E}$.  

In particular, it follows that 
\[
A\bousclass{\HFp } = \bigvee _{n\geq 0}\bousclass{T(n)} ,
\]
and hence
\[
A^{2}\bousclass{\HFp }=A \lbp \bigvee _{n\geq 0}\bousclass{T(n)} \rbp.
\]
We now prove that this is strictly larger that $\bousclass{\HFp }$,
using~\cite[Theorem 2.10]{rav-loc}.  Let
$J=(p^{i_{0}},v_{1}^{i_{1}},\dots ,v_{n}^{i_{n}},\dots )$ be an infinite
regular sequence in $\BP _{*}$.  Then we can form a spectrum $\bpj $
with $\bpj _{*}=\BP _{*}/J$ in various ways; Ravenel uses the
Bass-Sullivan construction.  By \cite[Corollary 2.14]{rav-loc}, $\bpj$
is a ring spectrum and hence is in $\DL$.  
Since $\bpj $ is built from $\BP $, we have
$\bpj\Smash A(n)=0$ for all $n$.  On the other hand, one can easily see
that $\bpj \Smash K(n)=0$ for all $n$, since a power of $v_{n}$ is
invariant modulo $(p^{i_{0}},v_{1}^{i_{1}},\dots ,v_{n-1}^{i_{n-1}})$
and this power has to act both invertibly and nilpotently on
$K(n)_{*}\bpj$.  Hence $\bousclass{\bpj}\leq A \lbp \bigvee _{n\geq
0}\bousclass{T(n)} \rbp$ for all infinite regular sequences $J$.  On the
other hand, Theorem~2.10 of~\cite{rav-loc} implies that, for almost all
such infinite regular sequences $J$, we have
$\bousclass{\bpj}>\bousclass{\HFp }$.
\end{proof}

In light of these results, we would like to understand
$A^{2}\bousclass{\HFp }$.  Given a regular sequence $J$ as in the
proof of Proposition~\ref{prop-H-is-not-closed}, we can form a
spectrum $S/J$ by taking the sequential colimit of the partial
quotients $S/J_{n}$.  This spectrum may not be well-defined even up to
Bousfield class, though each $S/J_{n}$ is.  The obvious conjecture is
that $A^{2}\bousclass{\HFp }$ should be the wedge of the
$\bousclass{S/J}$ over all infinite regular sequences $J$ and all
representatives $S/J$.

\section{Bousfield classes without finite acyclics}
\label{sec-finite-locals}

We have been discussing Bousfield classes with finite acyclics; in
this section, we examine the rest of the Bousfield classes.  No
spectrum can have both a nonzero finite acyclic and a nonzero finite
local; we conjecture that every spectrum has one or the other.  In any
case, we pay some attention to spectra with finite locals, and we
discuss Brown-Comenetz duality and its relation to such spectra.  We
also show that a number of conjectures related to Bousfield classes
without finite acyclics are equivalent.

Brown-Comenetz duality~\cite{brown-comenetz} is the main source
of counterexamples in the theory of Bousfield classes.  
Given a spectrum $X$, we denote by $\bcdual{X}$ its
\emph{Brown-Comenetz dual}, obtained by applying Brown
representability to the cohomology theory $Y\mapsto \Hom(\pi
_{0}(X\Smash Y),\mathbf{Q}/\mathbf{Z}_{(p)})$.  Let $I$ denote the
Brown-Comenetz dual of the sphere.  Note that $\bcdual{X}$ is the same
as the function spectrum $F(X,I)$, and there is a natural map
$X\xrightarrow{}\doublebcdual{X}$ which is an isomorphism when the
homotopy groups of $X$ are finite.  Also note that $\bcdual{X}=0$ if
and only if $X=0$, since $\mathbf{Q}/\mathbf{Z}_{(p)}$ is an injective
cogenerator of the category of $p$-local abelian groups.  The spectrum
$I$ is the central example of this paper.

Recall the spectra $X(n)$ from~\cite[Section 3]{rav-loc}, which
interpolate between the Bousfield classes of $S=X(0)$ and
$\BP=X(\infty)$: 
\[
\bousclass{S} = \bousclass{X(0)} > \bousclass{X(1)} > \dots >
\bousclass{X(\infty)} = \bousclass{\BP}.
\]

Some of the basic
properties of $I$ are as follows.

\begin{lemma}\label{lem-properties-of-I}\ \\ \vspace{-2ex}
\begin{enumerate}
\item $I$ is in the localizing subcategory generated by $\HFp$\usc{}
hence $\bousclass{\HFp} \geq \bousclass{I}$.
\item $X(1)\Smash I=0$\usc{} hence $X(n)\Smash I=0$ for all $n
\geq 1$, and $\BP\Smash I=0$.
\item $\HFp \Smash I = 0$\usc{} hence $\bousclass{\HFp} >
\bousclass{I}$, and $I \Smash I = 0$.
\item $T(n)\Smash I=0$ for all $n$. 
\item $\bousclass{I} \Smash \bousclass{F(n)} = \bousclass{IF(n)} =
\bousclass{I}$ for all $n$.
\item The mod $p$ Moore spectrum $M(p)$ \ulp and every
finite-dimensional torsion spectrum\urp \ is $I$-local.  
\end{enumerate}
\end{lemma}

\begin{proof}
Part~(a) follows immediately from the fact that the homotopy of $I$ is
bounded-above and torsion, as in \cite[2.6]{rav-loc}.  Part~(b) follows
from~\cite[Lemma 3.2]{rav-loc}, where it is shown that $[X(1),M(p)]=0$.
Using the isomorphism $M(p)=F(\bcdual{M(p)},I)$ and adjointness, we find
that $\bcdual{(X(1)\Smash \bcdual{M(p)})}=0$, so that $X(1)\Smash
\bcdual{M(p)}=0$.  Since the homotopy groups of $I$ are torsion, one can
readily verify that $\bousclass{\bcdual{M(p)}}=\bousclass{I}$, so that
$X(1)\Smash I=0$.  Since $\bousclass{\BP} \geq \bousclass{\BP \Smash
\HFp} = \bousclass{\HFp}$, then part~(c) follows from (a) and (b).
Part~(d) follows from part~(a) and the well-known fact that $\HFp \Smash
T(n)=0$ (because a $v_{n}$-self map must have positive Adams
filtration).  Part~(e) follows from (d) and the Bousfield class
decomposition of Section~\ref{sec-finite-acyclics}.  It is proved
in~\cite[Corollary B.13]{hovey-strickland} that $M(p)$ is $I$-local,
using the isomorphism $M(p)=\doublebcdual{M(p)}$.  It follows
from~\cite[Theorem B.6]{hovey-strickland} that every finite-dimensional
(defined in \cite{hovey-strickland}) torsion spectrum is $I$-local.
\end{proof}

Another useful property of $I$ is that it detects when a spectrum has
a finite local.  We have already discussed spectra with a finite
acyclic; similarly, we say that a spectrum $E$ \emph{has a finite
local} if there is a nonzero finite spectrum $X$ which is $E$-local.
Note that no spectrum can have both a nonzero finite local and a
nonzero finite acyclic: if $M$ is a finite $E$-local and $W$ is a
finite $E$-acyclic, then $M \Smash W$ is both local and acyclic, and
nonzero if both $M$ and $W$ are.  In ~\cite[Lemma
3.7]{hovey-chrom-split}, the first author shows that if $E$ has a
finite local, then every finite torsion spectrum is $E$-local.  This
was extended in~\cite[Theorem B.6]{hovey-strickland} to all
finite-dimensional torsion spectra.

\begin{proposition}\label{prop-I-detects-finite-locals}
The following are equivalent for a spectrum $E$\uc 
\begin{enumerate}
\item $M(p)$ is $E$-local. 
\item $E$ has a finite local. 
\item $aE\Smash I=0$. 
\item $\bousclass{E}\geq \bousclass{I}$. 
\end{enumerate}
\end{proposition}

\begin{proof}
We have already noted that (a) and (b) are equivalent.  To see that (c)
and (d) are equivalent, note that $aE\Smash I=0$ if and only if
$a\bousclass{E}\leq a\bousclass{I}$.  This holds if and only if
$\bousclass{E}\geq \bousclass{I}$.  

Since $M(p)$ is $I$-local, it follows that (d)$\Rightarrow$(a).  To see
that (a)$\Rightarrow$(c), suppose that $M(p)$ is $E$-local.  Then
$[aE,M(p)]_{*}=0$.  Using the isomorphism $M(p)=\doublebcdual{M(p)}$ and
adjointess, we find that $\bcdual{(aE\Smash \bcdual{M(p)})}=0$.  Thus
$aE\Smash \bcdual{M(p)}=0$.  We have already seen in the proof of
Lemma~\ref{lem-properties-of-I} that
$\bousclass{\bcdual{M(p)}}=\bousclass{I}$,
completing the proof.  
\end{proof}

Note that this proposition implies for example that every dissonant
spectrum is $I$-acyclic, since finite spectra are harmonic. 

\begin{example}\label{eg-meet-infinite-join}
Since finite spectra are harmonic, $\bousclass{I} < \bigvee
_{n}\bousclass{K(n)}$. In particular, $\bousclass{I}\meet \bigvee
_{n}\bousclass{K(n)}=\bousclass{I}$.  But for each $n$,
$\bousclass{I}\meet \bousclass{K(n)}=\bousclass{0}$, since
$\bousclass{K(n)}$ is minimal and $K(n)\Smash I=0$.  Thus the meet
does not distribute over infinite joins in the Bousfield lattice.
\end{example}


We now consider three conjectures, which we will prove are equivalent.
Note that for $X$ finite, $X \Smash I = F(\swdual{X},I) =
\bcdual{\swdual{X}}$, where $\swdual{X}$ is the Spanier-Whitehead dual
of $X$.  In particular, $X \Smash I \neq 0$ for every finite $X$.  This,
combined with Lemma~\ref{lem-properties-of-I}, suggests the following
conjecture, first made in~\cite[Appendix B]{hovey-strickland}.

\begin{conjecture}\label{conj-nothing-detects-I}
If $E\Smash I\neq 0$, then $\bousclass{E}\geq \bousclass{F(n)}$ for some
$n$.  
\end{conjecture}

Note that the converse to Conjecture~\ref{conj-nothing-detects-I}
is immediate from part~(d) of Lemma~\ref{lem-properties-of-I}.  

The following conjecture appeared in~\cite[Conjecture
3.10]{hovey-chrom-split}.

\begin{conjecture}[The Dichotomy Conjecture]\label{conj-dichotomy} 
Every spectrum has either a finite local or a finite acyclic.  
\end{conjecture}

It was pointed out in~\cite{hovey-coh-Bsfld} that the Dichotomy
Conjecture is equivalent to the following conjecture.

\begin{conjecture}\label{conj-everything-bigger-than-I}
If $E$ has no finite acyclics, then $\bousclass{E}\geq \bousclass{I}$.  
\end{conjecture}

The converse to Conjecture~\ref{conj-everything-bigger-than-I} follows
from Lemma~\ref{lem-properties-of-I}(e).

\begin{theorem}\label{thm-equivalence-of-conjectures}
The following are equivalent\uc
\begin{enumerate}
\item Conjecture~\ref{conj-nothing-detects-I}.
\item The Dichotomy Conjecture~\ref{conj-dichotomy}.
\item Conjecture~\ref{conj-everything-bigger-than-I}.
\end{enumerate}
\end{theorem}

\begin{proof}
We will prove that (a)$\Leftrightarrow$(b) and (b)$\Leftrightarrow$(c).

To see that (a)$\Rightarrow$(b), suppose that $E$ has no finite locals.
Then $aE\Smash I\neq 0$, by
Proposition~\ref{prop-I-detects-finite-locals}.  Hence, by part~(a),
$\bousclass{aE}\geq \bousclass{F(n)}$ for some $n$.  It follows that
$\bousclass{E}\leq \bousclass{a F(n)} = \bousclass{L_{n-1}^{f}S}$, and
so $E$ has a finite acyclic.

To see that (b)$\Rightarrow$(a), suppose $E\Smash I\neq 0$.  Then
$a^{2}E\Smash I\neq 0$, so $aE$ has no finite locals, again using
Proposition~\ref{prop-I-detects-finite-locals}.  Hence $aE$ must have a
finite acyclic, by part~(b), and so $\bousclass{aE}\leq
\bousclass{L_{n-1}^{f}S}$ for some $n$.  It follows that
$\bousclass{E}\geq \bousclass{F(n)}$ for some $n$.

Proposition~\ref{prop-I-detects-finite-locals} shows that
(b)$\Rightarrow$(c).  To see that (c)$\Rightarrow$(b), suppose that
$E$ has no finite acyclics.  Then part~(c) implies $\bousclass{E}\geq
\bousclass{I}$.  Since $M(p)$ is $I$-local by part~(e) of
Lemma~\ref{lem-properties-of-I}, it is also $E$-local.
\end{proof}

The Dichotomy Conjecture has a few interesting consequences.  The most
obvious one is that it implies that $\bousclass{I}$ is minimal.  

\begin{lemma}\label{lem-I-is-minimal}
If $E$ is a nontrivial spectrum with $\bousclass{E}<\bousclass{I}$,
then $E$ has no finite locals or finite acyclics.  Hence, if the
Dichotomy Conjecture holds, there are no such $E$, and $\bousclass{I}$
is a minimal Bousfield class.
\end{lemma}

\begin{proof}
Proposition~\ref{prop-I-detects-finite-locals} implies that $E$ has
no finite locals.  Since $\bousclass{E}<\bousclass{I}$, $E\Smash
L_{n-1}^{f}S=0$ for all $n$.  Thus $E$ can have no finite acyclics
either, by the Bousfield class decomposition of
Section~\ref{sec-finite-acyclics}.  
\end{proof}

The Dichotomy Conjecture also gives us a partial classification of
complemented Bousfield classes, when combined with the following
lemma.

\begin{lemma}\label{lem-complemented-detects-I}\ \\ \vspace{-2ex}
\begin{enumerate}
\item  Suppose that $K$ is a field spectrum.  Then for any $E$,
either $\bousclass{E}\geq \bousclass{K}$ or $\bousclass{aE}\geq
\bousclass{K}$.
\item At least one of $E$ and $aE$ has a finite local.  
\item If $E$ is complemented and has a finite local, then
$E\Smash I\neq 0$.
\end{enumerate}
\end{lemma}

\begin{proof}
(a): If $E\Smash K\neq 0$, then $\bousclass{E}\geq
\bousclass{E\Smash K}=\bousclass{K}$, since $E\Smash K$ is a nontrivial
wedge of suspensions of $K$.  If $E\Smash K=0$, then $\bousclass{K}\leq
\bousclass{aE}$ by definition of $\bousclass{aE}$.  

(b): Apply part~(a) to $\HFp $.  

(c): By Proposition~\ref{prop-I-detects-finite-locals}, since $E$
has a finite local, $aE\Smash I=0$.  Since $E$ is complemented, $aE$
must be its complement and $E\vee aE$ must detect every spectrum.  Thus
$E\Smash I\neq 0$.
\end{proof}


\begin{corollary}\label{cor-complemented}
None of the following spectra is complemented\uc{} $X(n)$, $\BP$,
$\HFp$, $\bigvee_{n}K(n)$, $\bigvee_{n}T(n)$, and $I$.  Furthermore,
if the Dichotomy Conjecture holds and $E$ is complemented, then either
$\bousclass{E}\geq \bousclass{F(n)}$ for some $n$ or
$\bousclass{E}\leq \bousclass{L_{n-1}^{f}S}$ for some $n$.
\end{corollary}

\begin{proof}
This follows from Lemma~\ref{lem-properties-of-I}.
\end{proof}

We have already seen that $K(n)$ is complemented for all $n$.  Hence
Corollary~\ref{cor-complemented} shows that $\BA$ is not closed under
infinite joins.  

By Proposition~\ref{prop-finite-acyclics-are-complemented},
Conjecture~\ref{conj-A(n)-minimal} implies the converse to the second
half of the corollary: every $E$ with $\bousclass{E}\geq
\bousclass{F(n)}$ or $\bousclass{E}\leq \bousclass{L_{n-1}^{f}S}$ is
complemented.  We can restate this as the following corollary.

\begin{corollary}\label{cor-BA}
Suppose both the Dichotomy Conjecture and
Conjecture~\ref{conj-A(n)-minimal} hold.  Then the
atoms of $\BA$ are $\bousclass{K(n)}$ and, for $n\geq 2$,
$\bousclass{A(n)}$.  Every element of $\BA$ can either be written as a
finite join of atoms or the complement of a finite join of atoms, in a
unique way.  In particular, $\BA$ is isomorphic to the Boolean algebra
of finite and cofinite subsets of a countable set.
\end{corollary}

\section{Strange Bousfield classes}\label{sec-strange}

In this section, we investigate some strange Bousfield classes.  We
start with the following problem.  As above, we write $\bcdual{X}$
for the Brown-Comenetz dual of $X$.

\begin{problem}\label{problem-strange}
Classify the strange Bousfield classes.  For instance, is every
strange spectrum Bousfield equivalent to $\bcdual{A}$ for some
connective $A$ with finitely generated homotopy groups?  Or to
$\bcdual{R}$ for some connective ring spectrum $R$?
\end{problem}

Note that $\bousclass{\bcdual{A}}\leq \bousclass{\HFp }$ for any
connective spectrum $A$ with finitely generated homotopy groups, since
then $\bcdual{A}$ will have homotopy groups bounded-above and torsion,
so will be in the localizing subcategory generated by $\HFp$.

While the set of strange Bousfield classes may be more complicated
than the guesses given in Problem~\ref{problem-strange}, these guesses
at least give us a starting place for the study of strange Bousfield
classes.  We find that when $A$ is as above, $\bcdual{A}$ is very much
like $I$.

\begin{lemma}\label{lem-R-finite-locals}
Suppose $A$ is a connective spectrum with finitely generated homotopy
groups.  Then the following are equivalent for a spectrum $E$.
\begin{enumerate}
\item $A\Smash M(p)$ is $E$-local. 
\item $A\Smash X$ is $E$-local for some finite torsion spectrum
$X$.  
\item $aE\Smash \bcdual{A}=0$. 
\item $\bousclass{E}\geq \bousclass{\bcdual{A}}$.
\end{enumerate}
\end{lemma}

The proof of this lemma is very similar to that of
Proposition~\ref{prop-I-detects-finite-locals}.  We require $A$ to
have finitely generated homotopy groups so that $A\Smash
X=I^{2} {(A\Smash X)}$ for all finite torsion $X$.  We require
that $A$ be connective as well so that $\bousclass{\bcdual{A}}\leq
\bousclass{\HFp }$.  This guarantees that $\bcdual{A}\Smash T(n)=0$
for all $n$, and thus that $\bousclass{\bcdual{(A\Smash
X)}}=\bousclass{\bcdual{A}}$ for all finite $X$.  We leave the rest of
the proof to the reader.




Similarly, we have the following analogue of
Theorem~\ref{thm-equivalence-of-conjectures}.  

\begin{theorem}\label{thm-A-equiv-of-conjs}
Suppose $A$ is connective and has finitely generated homotopy groups.
Then the following are equivalent. 
\begin{enumerate}
\item If $E\Smash \bcdual{A}\neq 0$, then $\bousclass{E}\geq
\bousclass{A} \Smash \bousclass{F(n)}$ for some $n$.  
\item For every $E$, either $A\Smash M(p)$ is $E$-local, or
$E\Smash A\Smash F(n)=0$ for some $n$.  
\item If $E\Smash A$ has no finite acyclics, then
$\bousclass{E}\geq \bousclass{\bcdual{A}}$.  
\end{enumerate}
\end{theorem}

Again we leave the proof to the reader.  The converses of parts~(a) and (c)
always hold, the key point being that $X\Smash \bcdual{X}$ is never
zero unless $X$ is.  Indeed, there is a map $X\Smash
\bcdual{X}\xrightarrow{}I$ adjoint to the identity map of
$\bcdual{X}$, and hence nontrivial.

We now examine some specific strange spectra.  We introduced the
spectra $X(n)$ in Section~\ref{sec-finite-locals}.

\begin{theorem}\label{thm-strange}
We have 
\[
\bousclass{I} =\bousclass{\bcdual{X(0)}} <\bousclass{\bcdual{X(1)}}<
\dots < \bousclass{\bcdual{X(n)}}< \bousclass{\bcdual{X(n+1)}}< \dots <
\bousclass{\bcdual{\BP}}< \bousclass{\HFp }.
\]
\end{theorem}

\begin{proof}
We first show that $\bousclass{\bcdual{X(n)}}\leq
\bousclass{\bcdual{X(n+1)}}$.  By Lemma~\ref{lem-R-finite-locals},
this is equivalent to showing that $X(n)\Smash M(p)$ is
$\bcdual{X(n+1)}$-local.  Because $X(n+1)\Smash
M(p)=\doublebcdual{(X(n+1)\Smash M(p))}$, we can use the same argument
as in the proof of Lemma~\ref{lem-properties-of-I}(e) to find that
$X(n+1)\Smash M(p)$ is $\bcdual{X(n+1)}$-local.  It therefore suffices
to show that $X(n)\Smash M(p)$ is in the colocalizing subcategory
generated by $X(n+1)\Smash M(p)$ (recall that a \emph{colocalizing
subcategory} is a thick subcategory closed under products).  We use
the $X(n+1)$-based Adams tower.  That is, we let $\overline{X(n+1)}$
be the fiber of the unit map of $X(n+1)$, we let
$X_{s}=\overline{X(n+1)}^{\Smash s}\Smash X(n)\Smash M(p)$, and we let
$K_{s}=X(n+1)\Smash X_{s}$.  There are then cofiber sequences
\[
X_{s+1}\xrightarrow{}X_{s}\xrightarrow{}K_{s} \xrightarrow{}\Sigma X_{s+1},
\]
and the homotopy inverse limit $\text{holim} (X_{s})$ is trivial for
connectivity reasons.  We turn this around by letting $X^{s}$ be the
cofiber of the map $X_{s}\xrightarrow{}X_{0}=X(n)\Smash M(p)$.  Then we
have cofiber sequences 
\[
X^{s+1}\xrightarrow{}X^{s} \xrightarrow{}\Sigma K_{s}
\xrightarrow{}\Sigma X^{s+1},
\]
and the homotopy inverse limit of $X^{s}$ is $X(n)\Smash M(p)$.  It
therefore suffices to show that each $K_{s}$ is in the colocalizing
subcategory generated by $X(n+1)\Smash M(p)$.

To see this, we use~\cite[Proposition 2.3]{devinatz-hopkins-smith},
which shows that $X(n+1)_{*}X(k)$ is a free module over $X(n+1)_{*}$ for
$k\leq n+1$.  It follows that $X(n+1)\Smash X(n+1)$ and $X(n+1)\Smash
X(n)$ are wedges of suspensions of $X(n+1)$.  Then one can easily check
that $K_{s}$ is a wedge of suspensions of $X(n+1)\Smash M(p)$, and since
everything is connective and locally finite, this wedge is also a
product.  Hence $K_{s}$ is in the colocalizing subcategory generated by
$X(n+1)\Smash M(p)$, and so $X(n)\Smash M(p)$ is as well.

A similar proof, using the fact that $\BP_{*}X(n)$ is a free
$\BP_{*}$-module, shows that $X(n)\Smash M(p)$ is in the colocalizing
subcategory generated by $\BP \Smash M(p)$.  Thus we have
$\bousclass{\bcdual{X(n)}}\leq \bousclass{\bcdual{\BP}}$.  We have
already seen that $\bousclass{\bcdual{X}}\leq \bousclass{\HFp }$ for any
connective $X$ with finitely generated homotopy groups.

It remains to show that all of the inequalities above are strict.  For
this we recall the method used by Ravenel in~\cite[Sections 2 and
3]{rav-loc}.  He shows that there are no maps from $X(n+1)$ to
$X(n)\Smash M(p)$, and that this is equivalent to the statement that 
\[
X(n+1)\Smash \bcdual{(X(n)\Smash M(p))}=0.
\]
We have already seen that
$\bousclass{\bcdual{X(n)}}=\bousclass{\bcdual{(X(n)\Smash M(p))}}$.
Hence $X(n+1)\Smash \bcdual{X(n)}=0$, but $X(n+1)\Smash \bcdual{X(n+1)}$ is
nonzero.  Thus $\bousclass{\bcdual{X(n)}}<\bousclass{\bcdual{X(n+1)}}$.
Similarly, Ravenel's proof that there are no maps from $\BP $ to
$X(n)\Smash M(p)$ shows that $\BP \Smash \bcdual{X(n)}=0$.  Since $\BP
\Smash \bcdual{\BP }$ is nonzero, this shows that
$\bousclass{\bcdual{X(n)}}<\bousclass{\bcdual{\BP }}$.  Finally, since
there are no maps from $\HFp $ to $\BP $, then $\HFp \Smash \bcdual{\BP }=0$,
and so $\bousclass{\bcdual{\BP }}<\bousclass{\HFp }$.  
\end{proof}

There are probably more strange Bousfield classes than the ones
described in Theorem~\ref{thm-strange}.  For example, Ravenel discusses
spectra $\bpj $ for infinite invariant regular sequences $J$ in $\BP
_{*}$ in ~\cite[Section 2]{rav-loc}.  We have already met these spectra
in the proof of Proposition~\ref{prop-H-is-not-closed}.  He shows that
$\bousclass{\BP }>\bousclass{\bpj }>\bousclass{\HFp }$ for $J\neq
(p,v_{1},\dots )$.  Presumably the Brown-Comenetz duals of these spectra
give other strange spectra.  In addition, at $p=2$, we have $\MSp$ as
well.  Ravenel sketched an argument to the first author once that
$\bousclass{\MSp}>\bousclass{\BP}$, and presumably one would also have
$\bousclass{\bcdual{\MSp}}<\bousclass{\bcdual{\BP }}$.

We do, however, make the following conjecture.  

\begin{conjecture}\label{conj-X(n)}
The spectra $X(n)$ and $X(n+1)$ are adjacent in the Bousfield lattice.
That is, if $\bousclass{E}> \bousclass{X(n+1)}$, then
$\bousclass{E}\geq \bousclass{X(n)}$.  
\end{conjecture}

Note that if Conjecture~\ref{conj-X(n)} holds, then
$a\bousclass{I}=\bousclass{X(1)}$.  Indeed, since $X(1)\Smash I=0$, we
have $a\bousclass{I}\geq \bousclass{X(1)}$.  Similarly, we have seen
above that $X(n+1)\Smash \bcdual{X(n)}=0$, so
$a\bousclass{\bcdual{X(n)}}\geq \bousclass{X(n+1)}$.  But $X(n)\Smash
\bcdual{X(n)}$ is nonzero, so we must have
$a\bousclass{\bcdual{X(n)}}=\bousclass{X(n+1)}$ if
Conjecture~\ref{conj-X(n)} holds.  Thus Conjecture~\ref{conj-X(n)} also
implies that $\bcdual{X(n)}$ and $\bcdual{X(n+1)}$ are adjacent in the
Bousfield lattice.  

Conjecture~\ref{conj-X(n)} also implies the following result. 

\begin{conjecture}\label{conj-xn-telescope}
$\bousclass{X(n)} \Smash \bousclass{T(k)}=\bousclass{T(k)}$ for all
$n$ and $k$.
\end{conjecture}

Hopkins has proved Conjecture~\ref{conj-xn-telescope}, but the authors
have not seen a proof.  To see that Conjecture~\ref{conj-X(n)} implies
Conjecture~\ref{conj-xn-telescope}, proceed by induction on $n$.  We
will only indicate the proof for $n=1$.  Conjecture~\ref{conj-X(n)}
implies that $\bousclass{X(1)} \vee \bousclass{T(k)}=\bousclass{X(1)}$.
By smashing with $T(k)$, we find that $\bousclass{T(k)}=\bousclass{X(1)}
\Smash \bousclass{T(k)}$, as required.

We mention that Hopkins has proved the following, though again the
authors do not know the proof. 

\begin{conjecture}\label{conj-cp-infty}
$\bousclass{S} = \bousclass{\mathbf{C}\!P^{\infty}}$.
\end{conjecture}

\section{Localizing and colocalizing subcategories}
\label{sec-localizing}

In this last section, we make a few remarks about general localizing
and colocalizing subcategories.  The outstanding question here is
whether every localizing subcategory is the collection of $E$-acyclics
for some $E$.  As pointed out by Neil Strickland, Ohkawa's
result~\cite{ohkawa} is relevant here.

Recall that a subcategory of the stable homotopy category is called
\emph{localizing} if it is thick and is closed under coproducts.  
The basic conjecture here is the following.

\begin{conjecture}\label{conj-localizing-subcats}
Every localizing subcategory is the collection of $E$-acyclics for some
$E$ \ulp{}and is therefore principal\urp{}.  
\end{conjecture}

There are several equivalent formulations of this conjecture.  First
we need some notation.  Given a spectrum $X$, let $\loc (X)$ denote
the localizing subcategory generated by $X$.

\begin{proposition}\label{prop-loc-subcat}
The following are equivalent.
\begin{enumerate}
\item Conjecture~\ref{conj-localizing-subcats} holds.
\item Every principal localizing subcategory $\loc (X)$ is the
collection of $E$-acyclics for some $E$.  
\item For each $X$, $\loc (X)$ is the collection of
$aX$-acyclics.
\item $\bousclass{X}\leq \bousclass{Y}$ if and only if $X\in
\loc (Y)$.
\end{enumerate}
\end{proposition}

\begin{proof}
It is clear that (a)$\Rightarrow$(b).  To see that (b)$\Rightarrow$(c),
suppose $\loc (X)$ is the $E$-acyclics for some $E$.  Then $E\Smash X=0$
so $\bousclass{E}\leq \bousclass{aX}$.  On the other hand, if $E\Smash
Z=0$, then $Z\in \loc (X)$, so $Z\Smash aX=0$.  Thus
$\bousclass{E}=\bousclass{aX}$, so $\loc (X)$ is also the collection of
$aX$-acyclics.  

To see that (c)$\Rightarrow$(d), note that $X\in \loc (Y)$ implies that
$\bousclass{X}\leq \bousclass{Y}$.  Conversely, if
$\bousclass{X}\leq \bousclass{Y}$, then $\bousclass{aX}\geq \bousclass{aY}$.
In particular, $X$ is an $aY$-acyclic.  Thus, from part~(c), we have
$X\in \loc (Y)$. 

It remains to show that (d)$\Rightarrow$(a).  We will first show
(d)$\Rightarrow$(c).  Indeed, suppose $Y$ is $aX$-acyclic.  Then
$\bousclass{Y}\leq \bousclass{aaX}=\bousclass{X}$.  By part~(d), we have
$Y\in \loc (X)$.  Hence $\loc (X)$ is the collection of $aX$-acyclics,
as required.  It is clear that (c)$\Rightarrow$(b), so it remains to
show that (b)$\Rightarrow$(a).  We will do so by showing that, given
part~(b), every localizing subcategory is principal.  Given a localizing
subcategory $\cat{C}$, there is only a set of Bousfield classes
represented by objects of $\cat{C}$ by~\cite{ohkawa}.  Since
(b)$\Rightarrow$(d), this means there is only a set of principal
localizing subcategories of $\cat{C}$.  Choose a representative for
each such principal localizing subcategory, and let $X$ be the wedge
of all of those representatives.  Then $\loc (X)=\cat{C}$, so
$\cat{C}$ is principal.
\end{proof}

Note that Conjecture~\ref{conj-localizing-subcats}, together with
Ohkawa's result, would imply that there is only a set of localizing
subcategories.  It would also imply that the cohomological
localizations studied in~\cite{hovey-coh-Bsfld} always exist, and are
in fact homological localizations.

We would like a similar understanding of colocalizing subcategories
(thick categories which are closed under products), but
such an understanding has eluded us.  The obvious conjecture is that
there is a one-to-one correspondence between localizing subcategories
and colocalizing subcategories, so that every colocalizing subcategory
would be the collection of $E$-locals for some $E$, given
Conjecture~\ref{conj-localizing-subcats}.  One could also ask whether
every colocalizing subcategory is principal.  We do not know the answer,
but we do have the following intriguing result.

Recall that a \emph{coideal} is a thick subcategory $\cat{C}$ with the
additional property that if $X\in \cat{C}$ and $Y$ is arbitrary, then
$F(Y,X)\in \cat{C}$.

\begin{proposition}\label{prop-coideal-I}
The colocalizing subcategory generated by $I$ is the entire stable
homotopy category, as is the coideal generated by $I$.
\end{proposition}

\begin{proof}
We use the results of~\cite{christensen-strickland}.  Recall that they
call a spectrum $X$ \emph{injective} if there are no phantom maps to it.
They show in~\cite[Proposition 3.9]{christensen-strickland} that
$\bcdual{X}$ is injective for all $X$.  They show in~\cite[Proposition
4.15]{christensen-strickland} that any $X$ fits into a cofiber sequence
$X\xrightarrow{}\doublebcdual{X}\xrightarrow{}K$, where $K$ is
injective.  It follows from~\cite[Lemma 4.14]{christensen-strickland}
that $K$ is a retract of $I^{2}K$.  Now, $\bcdual{Y}=F(Y,I)$ is in the
coideal generated by $I$ for any $Y$, so both $\doublebcdual{X}$ and $K$
are as well.  Hence $X$ is too.
\end{proof}



\providecommand{\bysame}{\leavevmode\hbox to3em{\hrulefill}\thinspace}

\end{document}